\documentclass[DIV=calc, paper=letter, fontsize=11pt]{scrartcl}
\usepackage[]{units}
\usepackage{lipsum} 
\usepackage[english]{babel} 
\usepackage[protrusion=true,expansion=true]{microtype} 
\usepackage{amsmath,amssymb,amsfonts,amsthm} 
\usepackage[svgnames, table]{xcolor} 
\usepackage[hang, small,labelfont=bf,up,textfont=it,up]{caption} 
\usepackage{booktabs} 
\usepackage{fix-cm}	 
\usepackage[]{units}

\usepackage{bm}
\usepackage{tabularx}

\usepackage[colorlinks=true, linkcolor=blue, citecolor=blue, filecolor=blue, runcolor=blue, urlcolor = blue]{hyperref}

\usepackage{fullpage}
\usepackage{sectsty} 

\usepackage{fancyhdr}
\usepackage{lastpage} 

\usepackage{nicefrac}

\usepackage{tikz}
\usetikzlibrary{calc,patterns,decorations.pathmorphing,decorations.markings}
\usepackage{pgfplots}

\usepackage{stmaryrd}
\usepackage{multicol}
\usepackage{float}

\usepackage{subcaption}

\usepackage{mathtools} 

\usepackage{amsmath,amssymb,amsfonts,amsthm}
\newtheorem*{remark}{Remark}

\newcounter{tempcolnum}
\makeatletter
\newcommand{\multicolinterrupt}[1]{
\setcounter{tempcolnum}{\col@number}
\end{multicols}
#1%
\begin{multicols}{\value{tempcolnum}}
}
\makeatother

\usepackage{graphicx}
\graphicspath{{figData/}}
\DeclareGraphicsExtensions{.png}

\usepackage{lettrine} 
\newcommand{\initial}[1]{ 
\lettrine[lines=3,lhang=0.3,nindent=0em]{
\color{DarkGoldenrod}
{\textsf{#1}}}{}}

\usepackage{caption}

\usepackage{graphicx}

\usepackage{titling} 

\newcommand{\HorRule}{\color{DarkGoldenrod} \rule{\linewidth}{1pt}} 

\pretitle{\vspace{-80pt} \begin{flushleft} \HorRule \fontsize{16}{16} \color{DarkRed} \selectfont} 

\title{Meshfree Methods on Manifolds for Hydrodynamic Flows on Curved Surfaces: A Generalized Moving Least-Squares (GMLS) Approach.} 

\posttitle{\par\end{flushleft}} 

\preauthor{\vspace{-5pt} \begin{flushleft}\fontsize{14}{14} \large  \color{DarkRed}} 
\author{B. J. Gross $^{1}$, N. Trask $^{3}$, P. Kuberry $^{3}$, and P. J. Atzberger$^{1,2}$ } 
\postauthor{\fontsize{12}{12} \small \color{Black} 
[$1$] Department of Mathematics, University of California Santa Barbara (UCSB); atzberg@gmail.com; http://atzberger.org/ [$2$] Department of Mechanical Engineering, University of California Santa Barbara (UCSB); [$3$] Sandia National Laboratories, Albuquerque, NM $^{**}$\\
\vskip 1.5em
\lfoot{\footnotesize * Work supported by DOE Grant ASCR PhILMS DE-SC0019246 and NSF Grant DMS-1616353. \\ **{\scriptsize Sandia National Laboratories is a multi-mission laboratory managed and operated by National
Technology and Engineering Solutions of Sandia, LLC., a wholly owned subsidiary of Honeywell International, Inc.,
for the U.S. Department of Energy’s National Nuclear Security Administration under contract DE-NA0003525. This
paper describes objective technical results and analysis. Any subjective views or opinions that might be expressed in
the paper do not necessarily represent the views of the U.S. Department of Energy or the United States Government. SAND2019-5752 J} }
\end{flushleft} \vspace{-18pt} \HorRule} 
\date{}  

\newcommand{\qq}[1]{#1}

\newcommand{\mb}[1]{\mathbf{#1}}
\newcommand{\bs}[1]{\boldsymbol{#1}}
\newcommand{\bsy}[1]{\boldsymbol{#1}}

\newcommand{\mDiv}{ {\mbox{div}} }
\newcommand{\mGrad}{ {\mbox{grad}} }
\newcommand{\mCurl}{ {\mbox{curl}} }

\newcommand{\h}{\mb{h}}

\DeclareGraphicsExtensions{.png}

\begin{document}

\maketitle 

\thispagestyle{fancy} 


\initial{W}\textbf{e utilize generalized moving least squares (GMLS) to develop meshfree techniques for discretizing hydrodynamic flow problems on manifolds.  We use exterior calculus to formulate incompressible hydrodynamic equations in the Stokesian regime and handle the divergence-free constraints via a generalized vector potential.  This provides less coordinate-centric descriptions and enables the development of efficient numerical methods and splitting schemes for the fourth-order governing equations in terms of a system of second-order elliptic operators.  Using a Hodge decomposition, we develop methods for manifolds having spherical topology.  We show the methods exhibit high-order convergence rates for solving hydrodynamic flows on curved surfaces.  The methods also provide general high-order approximations for the metric, curvature, and other geometric quantities of the manifold and associated exterior calculus operators.  The approaches also can be utilized to develop high-order solvers for other scalar-valued and vector-valued problems on manifolds.}

\setlength{\parindent}{5ex}

\section{Introduction} 
\label{sec_intro}

Many investigations in fluid mechanics pose challenges related to resolving hydrodynamic flows on curved surfaces or in confined geometries.  Examples include the transport of surfactants within bubbles and thin films~\cite{Hermans2015,Takagi2011,Kellay2017,Braun2012,SquiresManikantan2017}, protein drift-diffusion dynamics within lipid bilayer membranes and cell mechanics~\cite{Bassereau2011,Fuller2012,Saffman1975,Mogilner2018,Chou2010,Powers2002}, and colloidal aggregation within fluid interfaces~\cite{StebeCurvatureRodAssemblyPNAS2011,Lee2016,Ershov2013}.  Additional examples include stratified models in atmospheric and ocean science which employ shallow water equations within topologically spherical shells \cite{vallis2017atmospheric} and subsurface models governing the flow of groundwater through fractures in porous rock providing intricate geometries formed from the crack surfaces \cite{alboin2002modeling,chernyshenko2019unfitted,martin2005modeling,FrankGeological2007}.  \qq{For these problems the fluid mechanics can often be formulated in terms of effective fields on two dimensional surfaces.  In some cases these problems also can involve additional challenges of tracking an evolving geometry of the surface from the motion of the interface or even of tracking topological changes~\cite{Seifert1997,FrolovJanuary2015,SethianFoams2013}.}  \qq{We shall consider here primarily the problem of resolving hydrodynamic flows for surfaces of static shape.   Already in this case, challenges arise in formulating the hydrodynamic equations and numerical methods to take into account the contributions of the geometry.}

\qq{There has been a lot of interest in developing numerical methods to solve Partial Differential Equations (PDEs) on manifolds.}  Broadly categorized, these include Finite Element Methods (FEMs) ~\cite{DziukSurfFEM2007,BarrettFEMSurf2007,
DziukFEMSurfPDEs2013,DuSurfCH2011}, Level Set Methods (LSMs) and Phase Field Methods (PFMs)~\cite{bookOsherFedkiwLevelSets2003,
bookSethianLevelSets1996,Du2007,
OsherBertalmio2003,
SethianSmereka2003,
SethianFoams2013,
OsherBertalmio2003,Du2007}, Discrete Exterior Calculus Methods (DECMs)~\cite{DesbrunHirani2003,Hirani2003}, Finite Element Exterior Calculus Methods (FEECMs)~\cite{HolstExCalcFEM2017,Arnold2006,Arnold2010}, and other approaches~\cite{Zorin2005,
BertozziFourthOrderGeometries2006,
Mohammadi_Sphere_Kriging_LS_2019,
Stam2003,MacdonaldSurfPDEReactDiff2013,
LowengrubLeung2011}.  Each of these approaches have their strengths depending on the application addressed as well as having challenges.  FEMs offer specialized high-order methods with robust behaviors for broad problem classes with often rigorous guarantees of accuracy and stability when mesh quality factors for the geometry can be ensured~\cite{BrennerFEM2008}.  LSMs/PFMs provide an implicit representation of the geometry often more amenable to evolution and topological changes, but typically require sophisticated algorithms to track the interface, mitigate numerical diffusion, and recover quantities associated with the geometry and the scalar and vector fields on the surface~\cite{SethianSmereka2003,SethianFoams2013,OsherBertalmio2003,Du2007,OsherBertalmio2003,LowengrubLevelSetSurf2006,XuLevelSet2003}.  The DECMs/FEECMs provide discretizations with desirable qualities for mechanics allowing for derivation of methods that have conservation of mass, momentum, and vorticity~\cite{HiraniNavierStokes2016}.  By their design for preserving geometric structure, DECMs/FEECMs are currently applied primarily in fluid mechanics to inviscid flows.  While DECMs are elegant and very useful discretizations that have been applied successfully to many applications~\cite{HiraniHodgeStar2016,
DesbrunHirani2003,
HiraniNavierStokes2016,
DesbrunCrane2013}, for some scientific calculations they are low order, have limited convergence analysis~\cite{ErickSchulzDECConv2016,HaraniConvStudy2018}, or are restricted to specialized surface operations presenting some challenges for general physical modeling~\cite{Kotiuga2008,Bochev2006}.  In each of these methods, there is also a reliance upon a sufficiently high quality rectified or curvi-linear mesh or grid to locally represent the surface geometry or surface fields.  To complement these methods, we consider alternatives based on meshfree approaches for surface hydrodynamics and PDEs based on Generalized Moving Least Squares (GMLS) approximations~\cite{wendland2004scattered}.

\lfoot{} 

We develop GMLS approaches to approximate differential operators on manifolds where the shape is represented as a point set that samples the geometry.  We build on recent related work by 
Liang et al. who discretized the surface Laplace-Beltrami operator on manifolds~\cite{liang2013solving}.  We construct smooth continuous representations of the manifold by solving a collection of local least-squares problems over an approximating function space at each of the sample points to obtain local paramerizations.  This approach captures the geometry in a manner similar to~\cite{Shepard1968,Lancaster1981,Hoppe1992,Amenta2004}.  We approximate the surface scalar fields, vector fields, and differential operators by solving another collection of related local least-squares problems that make use of the geometric reconstructions.  In conjunction, these provide general methods for obtaining high-order approximations of the manifold shape, operators arising in differential geometry, and operators of differential equations.  We use exterior calculus for generalizing operations from vector calculus and techniques from mechanics to the manifold setting.  This provides a convenient way to formulate incompressible hydrodynamic equations for flows on curved surfaces and related GMLS approximations.  \qq{We also use these approaches to show in general how equations and related solvers can be formulated in terms of vector potentials facilitating development of other physical models with constraints and numerical methods.}

\qq{We also mention there are many existing meshfree approaches for solving PDEs.  These may be characterized broadly by the underlying discretizations.  This includes Radial Basis Functions (RBFs) \cite{buhmann2003radial}, Smooth Particle Hydrodynamics (SPH) \cite{gingold1977smoothed}, and the approaches of Generalized Finite Differences / Moving Least Squares / Reproducing Kernel Particle Method (GFD/MLS/RKPM) \cite{lancaster1981surfaces}.}  While the majority of meshfree literature has concerned solution of PDEs in $\mathbb{R}^d$, significant recent work has focused on the manifold setting~\cite{LowengrubLeung2011,
liang2013solving,
suchde2019,Shankar_RBF_Lagrangian_2018,
Mohammadi_Sphere_Kriging_LS_2019,
MacdonaldSurfPDEReactDiff2013,
ArroyoSurfMonge2019,PiretRBF_Orthog_Grad2012,FuselierWright2013}.  In the last decade, substantial work has been done to use RBFs to solve shallow-water equations on the sphere~\cite{flyer2007transport}.  \qq{The meshfree setting is attractive particularly for building semi-Lagrangian schemes of interest in atmosphere science and other applications~\cite{Shankar_RBF_Lagrangian_2018}. Significant work also has been done on RBF methods to obtain robust numerical methods for predictive simulations in~\cite{fornberg2011stabilization,flyer2012guide,fries2008convergence} and for solving PDEs on manifolds without the need for local surface reconstructions in~\cite{ShankarRBF_FD_RD_2015,ShankarRBF_LOI_Aug_SurfPDE2018}.  Recent work on RBF-FD also includes methods for reaction-diffusion equations on surfaces and other PDEs~\cite{ShankarRBF_FD_RD_2015,LehtoWrightRBF_Compact2017,ShankarRBF_LOI_Aug_SurfPDE2018} and related approaches in~\cite{PiretRBF_Orthog_Grad2012,PetrasClosestPointRBF2018,Fornberg_RBF_Flat_Stable2008}.}

\qq{SPH approaches have also been introduce that offer attractive structure-preserving properties, particularly in conserving invariants of Lagrangian transport.  However, in general it is not possible for SPH to simultaneously obtain conservation principles and a consistent discretization~\cite{trask2015scalable}.}  \qq{The MLS/RKPM/GFD approaches provide a compelling alternative by addressing accuracy issues through the explicit construction of approximations with polynomial reproduction properties and an accompanying rigorous approximation theory~\cite{wendland2004scattered,schaback2016error}.  However, it should be noted in many cases stability theory currently is still lacking.  There have been several examples of successful discretizations for scalar surface PDEs in~\cite{sokolov2017numerical,suchde2018meshfree}.}  

\qq{In Generalized Moving Least Squares (GMLS) this approach is extended to enable the recovery of arbitrary linear bounded target functionals from scattered data~\cite{wendland2004scattered,mirzaei2012generalized}.  For transport and flow problems in $\mathbb{R}^d$, compatible GMLS methods have been developed in~\cite{trask2018compatible} which parallel the stability of compatible spatial discretization~\cite{arnold2007compatible}.  In the Euclidean setting, this has allowed for stable GMLS discretizations of Darcy flow in $\mathbb{R}^d$\cite{trask2017high}, Stokes flow in $\mathbb{R}^d$ \cite{trask2018compatible}, and fluid-structure interactions occurring in suspension flow \cite{hu2019spatially}.  In the recent work~\cite{traskcompatible}, is has been shown that the scheme developed by Liang et al.~\cite{liang2013solving} to discretize the Laplace-Beltrami operator on manifolds admits an interpretation as a GMLS approximation.  This unification enables extensions of the compatible staggered approach for Darcy in $\mathbb{R}^d$ \cite{trask2017high} to the manifold setting~\cite{traskcompatible}.}

\qq{We develop here related methods for discretizing the diverse collection of exterior calculus operators to obtain high-order solutions to PDEs on surfaces.}  We focus particularly on the case of developing methods for hydrodynamic flows on curved surfaces.  We introduce background on the GMLS approximation approach in Section~\ref{sec:gmls}.  In Section~\ref{sec:geo_reconstruct}, we discuss how to use GMLS to reconstruct locally the manifold geometry from a point set representation, approximate quantities from differential geometry, and approximate operators that generalize vector calculus to the manifold setting.  \qq{In Section~\ref{sec:hydro_formulation},
we show how exterior calculus approaches can be used to formulate equations for hydrodynamic flow on surfaces in a few different ways which facilitates development of alternative solvers.}  We discuss our numerical solvers for incompressible hydrodynamic flows in Section~\ref{sec:numerical_solver}.  Finally, in Section~\ref{sec:results}, we conclude with results discussing our investigations of the accuracy of the GMLS methods.  \qq{In particular, we study convergence of the approximations for the operators on the manifold and the precision of our solvers for hydrodynamic flows on curved surfaces.  Many of our methods can be adapted readily for approximating other scalar-valued and vector-valued PDEs on manifolds.}

\section{Generalized Moving Least Squares (GMLS)}
\label{sec_gmls}
\label{sec:gmls}
\qq{The method of Generalized Moving Least Squares (GMLS) is a non-parametric functional regression technique for constructing approximations by solving a collection of local least-squares problems based on scattered data samples of the action of a target operator~\cite{wendland2004scattered,mirzaei2012generalized,mirzaei2012generalized,Mirzaei2016}.  These local problems are formulated by specifying a finite collection of functionals that probe features of the action of the target operator.}

\qq{More specifically, consider a Banach space $\mathbb{V}$ and function $u \in \mathbb{V}$.  We assume that $u$ is characterized by a scattered collection of sampling functionals $\Lambda(u) := \left\{\lambda_j(u)\right\}_{j=1}^N \subset \mathbb{V}^*$, where $\mathbb{V}^*$ is the dual of $\mathbb{V}$.  Here, we shall primarily use sampling functionals that are point evaluations $\lambda_i(u) = \delta_{x_i} [u] = u(x_i)$. We denote the collection of sample points as $\mathbb{X}_h := \left\{\mathbf{x}_j\right\}_{j=1}^N$, where $h$ indicates the spatial resolution. We assume $\mathbb{X}_h \subset \Omega \subset \mathbb{R}^d$ for a compactly supported domain $\Omega$.  We characterize the distribution of points by 
\begin{equation}
\label{eqn:def_q_uniform}
h_{\mathbb{X},\Omega} = \underset{\mathbf{x} \in \Omega}{\sup}\underset{1 \leq j \leq N}{\min}||\mathbf{x} - \mathbf{x}_j||_2, 
\hspace{1cm}
q_{\mathbb{X}} = \frac12 \underset{i \neq j}{\min} ||\mb{x}_i - \mb{x}_j||_2,
\hspace{1cm}
q_{\mathbb{X}} \leq h_{\mathbb{X},\Omega} \leq c_{qu} q_{\mathbb{X}}.
\end{equation}
The $||\cdot||_2$ is the Euclidean norm, 
$h_{\mathbb{X},\Omega}$ is the \textit{fill distance},
$q_{\mathbb{X}}$ is the \textit{separation distance} of $\mathbb{X}_h$. The point set is called \textit{quasi-uniform} if there there exists $c_{qu} > 0$ in the last expression of equation~\ref{eqn:def_q_uniform}.  We shall assume $\mathbb{X}_h$ is quasi-uniform throughout, which is important in proving results about existence, convergence, and accuracy of GMLS~\cite{wendland2004scattered,mirzaei2012generalized}.}

\qq{Consider a target linear functional $\tau_{\hat{\mb{x}}}$ at location $\hat{\mb{x}}$. For example, the point-evaluation of a differential operator $\tau_{\hat{\mb{x}}} = D^\alpha u(\hat{\mb{x}})$ with
$\alpha$ the multi-index.  To approximate such operators, we first solve using the samples   a collection of local weighted $\ell_2$-optimization problems over a finite dimensional subspace $\mathbb{V}_h \subset \mathbb{V}$.  In particular, we solve for $p^* \in \mathbb{V}_h$ with
\begin{eqnarray}
\label{eqn:opt_gmls_reconstruct}
\label{eqn:weight_func}
    p^* = \underset{{q \in \mathbf{V}_h}}{\mbox{argmin}} \sum_{j=1}^N \left( \lambda_j(u) -\lambda_j(q) \right)^2 \omega(\lambda_j,\tau_{\tilde{\mb{x}}}),
\hspace{1cm}
\omega(\lambda_j,\tau_{\tilde{\mb{x}}}) = \Phi(||\mb{x}_j - \tilde{\mb{x}}||_2).
\end{eqnarray}}
\qq{The $\omega$ is a compactly supported positive function correlating information at the sample location $\mb{x}_j$ and the target location $\tilde{\mb{x}}$.  Throughout, we take $\Phi$ to be radially symmetric with  $\Phi(r) = \left(1 - {r}/{\epsilon}\right)^{\bar{p}}_+$, where $(z)_+ = \max\{z,0\}$ and $\bar{p} > 0$ with the $\epsilon$ controlling the shape and support of $\omega$.}

\qq{For the basis $\mathbb{V}_h = \mbox{span}\{\phi_1,...,\phi_{\dim(\mathbb{V}_h)} \}$, we denote by $\mathbb{P}(x)$ the vector whose $i^{th}$-entry is $[\mathbb{P}(x)]_i = \phi_i(x)$.}  \qq{The solution to equation \ref{eqn:opt_gmls_reconstruct} can be represented using a coefficient vector $\mathbf{a}(u)$ to
express the GMLS approximation of $\tau_{\tilde{\mb{x}}}$ as
\begin{equation}
\label{eqn:gmls_approx_tau1}
 p^* = P(x)^\intercal \mathbf{a}(u),
\hspace{1cm}
\tau^h_{\tilde{\mb{x}}}(u) = \tau_{\tilde{\mb{x}}}(\mathbb{P})^\intercal \mathbf{a}(u).
\end{equation}
}
\qq{Assuming that the collection of sampling functionals $\Lambda$ is unisolvent for $\mathbb{V}_h$, the GMLS estimate of $\tau_{\tilde{\mb{x}}}$ in equation \ref{eqn:gmls_approx_tau1} can be expressed as
\begin{equation}
\label{eqn:gmls_approx_tau2}
\tau^h_{\tilde{\mb{x}}}(\phi) = \tau_{\tilde{\mb{x}}}(\mathbf{P})^\intercal \left(\Lambda(\mathbf{P})^\intercal \mathbf{W} \Lambda(\mathbf{P})\right)^{-1} \Lambda(\mathbf{P})^\intercal \mathbf{W} \Lambda(u).
\end{equation}
}
\qq{The $\Lambda$ is called \textit{unisolvent} over $\mathbb{V}_h$, if  any element of $\mathbb{V}_h$ is uniquely determined by the collection of sampling functionals $\lambda_j$, here by the points in the support of $\omega$ \cite{wendland2004scattered}.}

\qq{We summarize the GMLS approximation approach in Figure~\ref{fig:gmls_overview}.   We use the following notation throughout
\begin{itemize}
  \item $\tau_{\tilde{\mb{x}}}(\mathbf{P}) \in \mathbb{R}^{\dim(V_h)}$ denotes the vector with components consisting of the target functional applied to each of the basis functions $\phi_k$.
  \item $\mathbf{W} \in \mathbb{R}^{N \times N}$ denotes the diagonal weight matrix with entries $\left\{\omega(\lambda_j,\tau_{\tilde{\mb{x}}})\right\}_{j=1}^N$.
  \item $\Lambda(\mathbf{P}) \in \mathbb{R}^{N \times \dim(V_h)}$ denotes the rectangular matrix whose $(j,k)$-entry is $\lambda_j(\phi_k)$ corresponds to the application of the $j^{th}$ sampling functional $\lambda_j$ applied to the $k^{th}$ basis function $\phi_k$.
  \item $\Lambda(u) \in \mathbb{R}^N$ denotes the vector consisting of entries $\{\lambda_j(u)\}_{j=1}^N$ corresponding to the $N$ sampling functionals $\lambda_j$ applied to the function $u$.
\end{itemize}
}

\qq{We remark that an advantage of GMLS over traditional least-squares approaches is that to build up approximations it only requires information locally at nearby points.  Algorithmically, the main expense in GMLS is in inverting over the base points $\tilde{\mb{x}}$ many separate small systems of dense normal equations given by equation~\ref{eqn:gmls_approx_tau2}.}  The GMLS approach is very well-suited to hardware acceleration and parallelization using packages such as the recent Compadre toolkit~\cite{Compadre2019}.  

We shall consider here primarily the case when the target functional $\tau$ is selected to approximate point evaluations of either the function as in regression or of differential operators acting on manifolds.  These approximations and estimates have relation to~\cite{Nayroles1992,mirzaei2012generalized,wendland2004scattered,Mirzaei2016}.  \qq{For partial derivatives $\mathcal{D}^{\bsy{\alpha}}$ in $\mathbb{R}^d$ with multi-index $\bsy{\alpha}$, Mirzaei~\cite{mirzaei2012generalized} proved the estimate }
\begin{eqnarray}
\label{eqn:GMLS_deriv_bound}
\|\mathcal{D}^{\bsy{\alpha}} u - \mathcal{D}^{\bsy{\alpha}}p^*\|_2 \leq C h^{m + 1 - |\alpha|}.
\end{eqnarray}
These reconstructions use $m^{th}-$order polynomials. We extend such approximations to the manifold setting to handle non-linear target functionals due to geometry-dependent terms.  \qq{While Mirzaei's analysis was not developed for the non-linear setting, we find empirically convergence rates manifest in our approximations similar to equation~\ref{eqn:GMLS_deriv_bound}.}

\begin{remark}
\qq{We consider throughout quasi-uniform point sets of two-dimensional compact manifolds embedded in $\mathbb{R}^3$. It can be shown readily~\cite{wendland2004scattered} that for the quasi-uniform Euclidean setting in $\mathbb{R}^2$ that there exists constants $c_1,c_2 > 0$ such that
$c_1 h_{\mathbb{X},\Omega} \leq \frac{1}{\sqrt{n}} \leq c_2 h_{\mathbb{X},\Omega}$, and thus the fill distance scales as $h \sim 1/\sqrt{n}$, where $n$ is the number of points. We shall use the notation $\bar{h}^{-1} := \sqrt{n}$ to characterize the refinement level of the point set.}
\end{remark}

\section{Geometric Reconstructions from the Point Set of the Manifold}
\label{sec:geo_reconstruct}

\begin{figure}[H]
\centerline{\includegraphics[width=1.00\columnwidth]{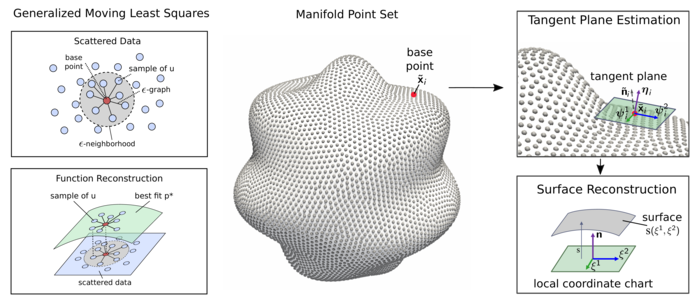}}
\caption{\qq{GMLS Approximation of Operators and Surface Reconstructions.  A target functional 
$\tau_{\tilde{\mb{x}}}[u]$ is approximated using 
data within an $\epsilon$-neighborhood around the 
base point $\tilde{\mb{x}}$ \textit{(left-top)}.  
For values of $u$ the best fitting function 
$p^* \in \mathbb{V}_h$ is identified using the sampling functionals 
$\{\lambda_j\}$ for computing
$\tau_{\tilde{\mb{x}}}^h[u] = \tau_{\tilde{\mb{x}}}[p^*]$~\textit{(left-bottom)}.  
For geometric reconstructions~\textit{(center)}, a Principle Component Analysis (PCA) is used to find local parameterization of the surface of form $(\xi^1,\xi^2,s(\xi^1,\xi^2))$~\textit{(top-right)}.  The $s(\xi^1,\xi^2)$ 
and its derivatives are approximated by GMLS to obtain general geometric quantities of the manifold and approximate operators~\textit{(bottom-right)}.}}
\label{fig:mesh_reconstruct}	
\label{fig:gmls_overview}	
\end{figure}

\qq{To formulate GMLS problems on manifolds, we must develop estimates of the metric tensor and other geometric quantities associated with the shape of the manifold.  The metric tensor and geometric quantities are first extracted from the point cloud sampling of the manifold and then used to approximate the differential operators on the surface.}

\qq{Consider a smooth manifold $\mathcal{M} \subset \mathbb{R}^d$ and assume a quasi-uniform point cloud sampling $\mathbb{X}_h \subset \mathcal{M}$.}  At each point $\mb{x}_i \in \mathbb{X}_h$, we shall construct an approximation to the tangent space $T_{\mathbf{x}_i}$ \cite{wendland2004scattered,liang2013solving}.  \qq{At location $\mb{x}_i$, we use Principal Component Analysis (PCA) based on nearby neighbor points $\mb{x}_j$ such that $j \in \mathcal{N}_i$.  The $\mathcal{N}_i = \mathcal{N}_\epsilon(\mb{x}_i) = \{\mb{x}_k \,|\, \|\mb{x}_k - \mb{x}_i\|_2 < \epsilon \} = \mathbb{X}_h \bigcap B_\epsilon(\mb{x}_i)$ is the $\epsilon-$ball of samples around $\mb{x}_i$.}  \qq{To center the sample points for use in PCA, we define the centering point}
\begin{equation}
    \bar{\mb{x}}_i = \frac{1}{|\mathcal{N}_i|} \underset{j \in \mathcal{N}_i}{\sum} \mb{x}_j.
\end{equation}
\qq{While in general $\bar{\mb{x}}_i \neq \mb{x}_i$, in practice these are typically close. We also refer to $\mathcal{N}_i = \mathcal{N}_\epsilon(\mb{x}_i)$ as the patch of points at $\mb{x}_i$.  For PCA we use for the empirical estimate of the covariance at $\mb{x}_i$}
\begin{equation}
    \mathcal{C} = \mbox{Cov}(\{\mb{x}_j\}) = \frac{1}{|\mathcal{N}_i|} \underset{j \in \mathcal{N}_i}{\sum} \left(\mb{x}_j - \bar{\mb{x}}_i\right)\left(\mb{x}_j - \bar{\mb{x}}_i\right)^\intercal.
\end{equation}
\qq{This provides a good estimate to the local geometry when $h_{\mathbb{X}}$ and $\epsilon$ are sufficiently small that the set of points $\mathcal{N}_\epsilon(\mb{x}_i)$ is nearly co-planar.}  We estimate the tangent space $T\mathcal{M}_{\mb{x}_i}$ of the manifold using the $(d-1)$-largest eigenvectors of $\mathcal{C}$.  These provide when $d = 3$ a basis for the tangent plane that we denote by $\bsy{\psi}^1_i$ and $\bsy{\psi}^2_i$ and normalize to have unit magnitude.  These also give the unit normal as $\bsy{\eta}_i = \bsy{\psi}^1_i \times \bsy{\psi}^2_i$.  \qq{We show the steps in the geometric reconstruction approach in Figure~\ref{fig:mesh_reconstruct}.}  

\begin{remark}
It is important to note that the PCA-approach can arbitrarily assign an orientation in the reconstruction of the tangent space.  This can have the undesirable property that neighboring patches have opposite orientations resulting in sign changes for some surface operators, such as the curl.  
In the general case, globally orienting the surface is a challenging NP-hard problem, as discussed in Wendland~\cite{wendland2004scattered}.  Many specialized algorithms have been proposed for this purpose which are efficient in practice, including front-marching and voronoi-based methods~\cite{DesbrunNormals2007,wendland2004scattered}.  We shall assume throughout that at each point $\mb{x}_i$ there is a reference normal $\tilde{\mb{n}}_i$ either determined in advance algorithmically or specified by the user.  We take in our PCA procedures that the normals $\bsy{\eta}_i$ are oriented with $\tilde{\mb{n}}_i^T \bsy{\eta}_i > 0$.
\end{remark}

We use this approach to define a local coordinate chart for the manifold in the vicinity of the base point $\tilde{\mb{x}} = \mb{x}_i$.  \qq{For this purpose, we take as the origin the base point $\mb{x}_i$ and use the tangent plane bases $\bsy{\psi}_i^1, \bsy{\psi}_i^2$ and normal $\bsy{\eta}_i$ obtained from the PCA procedure.}  We then define a local coordinate chart using the embedding map $\bsy{\sigma}$
\begin{eqnarray}
\label{eqn:chart}
\bsy{\sigma}(\xi^1,\xi^2;q) = \mb{x}_i + \xi^1 \bsy{\psi}_i^1 + \xi^2 \bsy{\psi}_i^2 + s(\xi^1,\xi^2)\bsy{\eta}_i.
\end{eqnarray}
This provides a family of parameterizations in terms of local coordinates  $(\xi^1,\xi^2)$, defined by choice of a smooth function $s$.
 Without loss of generality we could always define the ambient space coordinates so that locally at a given base point $\tilde{\mb{x}}$ we have $\bsy{\sigma} = (\xi^1,\xi^2,s(\xi^1,\xi^2))$.  This can be interpreted as describing the surface as the graph of a function over the $(\xi^1,\xi^2)$-plane where $s$ is the height above the plane, see Figure~\ref{fig:mesh_reconstruct}.  This parameterization is known as the Monge-Gauge representation of the manifold surface \cite{Monge1807,Pressley2001}, and we will use GMLS to approximate derivatives of $\bsy{\sigma}$ through the following choices:
 
\begin{itemize}
\item We take for our sampling functionals $\Lambda = \{\lambda_j\}_{j=1}^N$ point evaluations $\lambda_j = \delta_{\mb{x}_i}$ at all points $\mb{x}_j$ in the $\epsilon$-ball neighborhood $\mathcal{N}_i$ of $\mathbf{x}_i$.
\item \qq{We use the target functional $\tau^{[\bsy{\alpha}]}$ as the point evaluation of the derivative $\mathcal{D}^{\bsy{\alpha}} \bsy{\sigma}$ at $\mathbf{x}_i$, where $\mathcal{D}^{\bsy{\alpha}}$ denotes the partial derivative of $\bsy{\sigma}$ in $\{\xi^c\}$ described by the multi-index $\bsy{\alpha}$ \cite{evans2010partial}.} 
\item We take for the reconstruction space the collection of $m_1^{th}$-order polynomials. 
\item We use for our weighting function the kernel in equation~\ref{eqn:weight_func} with support matching the parameter $\epsilon$ used for selecting neighbors in our reconstruction and for defining our $\epsilon$-graph on the point set.
\end{itemize}

We use these point estimates of the derivative of $\bsy{\sigma}$ to evaluate non-linear functionals of $\bsy{\sigma}$ characterizing the geometry of the manifold. Consider the metric tensor
\begin{eqnarray}
g_{ab} = \langle \bsy{\sigma}_{\xi^a}, \bsy{\sigma}_{\xi^b} \rangle_g.
\end{eqnarray}
The $\langle \mb{a}, \mb{b} \rangle_g$ corresponds to the usual Euclidean inner-product $\mb{a} \cdot \mb{b}$ when the vectors $\bsy{\sigma}_{\xi^c} = \partial \bsy{\sigma}/\partial \xi^c$ are expressed in the basis of the ambient embedding space. Other geometric quantities can be similarly calculated from this representation once estimates of $\mathcal{D}^{\bsy{\alpha}} \bsy{\sigma}$ are obtained.

\subsubsection{GMLS Approximation of Geometric Quantities}
\label{sec:geo_approx_quant}

\begin{figure}[H]  
\centerline{\includegraphics[width=0.8\columnwidth]{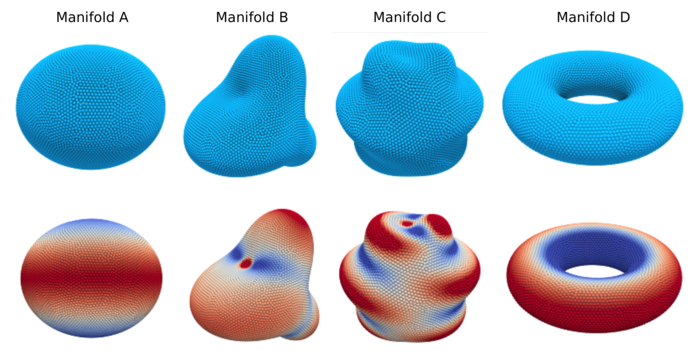}}
\caption{\qq{Point Set Representations of Manifolds~\textit{(bottom)}.  Manifold A is an ellipsoid defined by the equation ${x^2}/{a^2} + {y^2}/{b^2} + z^2 = s_0^2$ with $a = 1.2, b = 1.2, s_0^2 = 1$.  Manifold B is a radial manifold defined in spherical coordinates by $(\theta,\phi,r(\theta,\phi)$ where $r(\theta,\phi) = 1 + r_0\sin(3\phi)\cos(\theta)$ with $r_0 = 0.1$.  Manifold C is a radial manifold defined in spherical coordinates by $(\theta,\phi,r(\theta,\phi)$ where $r(\theta,\phi) = 1 + r_0\sin(7\phi)\cos(\theta)$ with $r_0 = 0.1$.  Manifold D is a torus defined by the equation $(s_1^2 - \sqrt{x^2 + y^2} )^2 + z^2 = s_2^2$ with $s_1^2 = 0.7, s_2^2 = 0.3$.  Estimated Gaussian curvature of the manifold~\textit{(bottom)}.
Each of the manifolds shown are represented by quasi-uniform point sets with approximately $n = 10^4$ samples.  For quasi-uniform sampling we expect the fill-distance $h$ to scale as $h \sim 1/\sqrt{n}$.  When reporting our results, we use throughout the notation $\bar{h}^{-1} = \sqrt{n}$.  We discuss further details of the point sampling of the manifolds in Appendix~\ref{appendix:sampling_res_manifolds}.}
}
\label{fig:manifolds}	
\label{fig:gauss_curv_surf}	
\end{figure}

\qq{We now utilize these approaches to estimate Gaussian curvature, as a representative geometric quantity of interest. In Appendix~\ref{appendix:monge_gauge} we provide detailed expressions for additional geometric quantities of interest which we will later need for approximating hydrodynamic flows on surfaces.} To demonstrate in practice the convergence behavior of our techniques as the fill-distance is refined, we consider the four example manifolds shown in Figure~\ref{fig:manifolds}.

We utilize the Weingarten map $\mb{W} = \mb{I}^{-1} \mb{II}$ to estimate the Gaussian curvature via the formula $K = \det(\mb{W})$ when using the GMLS estimate of $\bsy{\sigma}_{\xi^c}$ to calculate $\mb{I}$ and $\mb{II}$, see Appendix \ref{appendix:monge_gauge}.  \qq{We investigate the convergence of the estimated curvature for the manifolds A-D as the point sampling resolution increases in Table~\ref{table:conv_Gauss_Curv}.  We show the estimated curvature on the surface for each of the manifolds in Figure~\ref{fig:gauss_curv_surf}.}  We find our GMLS methods with $m = 6$ yields approximations having $5^{th}$-order accuracy.  While there is currently no convergence theory for our non-linear estimation procedure, the results for $k = 2$ for Gaussian Curvature are consistent with the suggestive predictions $m + 1 - k$ of equation~\ref{eqn:GMLS_deriv_bound}.

\begin{table}[H]
\label{K_Conv}
\begin{center}
\begin{tabular}{ccccccccccc}
\textbf{}  & \multicolumn{2}{c}{\textbf{Manifold A}}      & \multicolumn{2}{c}{\textbf{Manifold B}} & \multicolumn{2}{c}{\textbf{Manifold C}} & & & \multicolumn{2}{c}{\textbf{Manifold D}}      \\\cline{2-7}\cline{10-11}
\textbf{h} & \textbf{$\ell_2$-error} & \textbf{Rate} & \textbf{$\ell_2$-error} & \textbf{Rate} & \textbf{$\ell_2$-error} & \textbf{Rate} & & \textbf{h} & \textbf{$\ell_2$-error} & \textbf{Rate} \\
0.1 & 2.1351e-04 & - & 1.1575e-01 & -  & 1.2198e-01  & -&  & .08 &  5.5871e-02 & -  \\
0.05 & 3.0078e-06 & 6.07 & 1.6169e-02 & 2.84 & 4.7733e-03 & 4.67 & & .04 & 6.5739e-04 & 6.51 \\
0.025 & 5.3927e-08 & 5.77 & 8.3821e-04 & 4.26 & 1.6250e-04 & 4.88 & & .02 &  1.3418e-05 & 5.67 \\
0.0125 & 1.1994e-09 & 5.48 & 2.3571e-05 & 5.14 & 4.5204e-06 &  5.17 & & .01 & 3.1631e-07 & 5.37 \\
\end{tabular}
\end{center}
\caption{\qq{Convergence of GMLS Approximation of the Gaussian Curvature $K$.  The GMLS reconstruction of the manifold is used with polynomial order $m_1 = 6$.  Our GMLS methods involve operations with $k_1 = 2^{nd}$-order differentiation.  We find $~\sim 5^{th}$-order asymptotic convergence rate.  The target sampling distance $h$ is discussed in Appendix~\ref{appendix:sampling_res_manifolds}.  The estimated Gaussian curvature for each manifold is shown in Figure~\ref{fig:gauss_curv_surf}.}}
\label{table:conv_Gauss_Curv}
\end{table}

\subsection{Generalizing the Differential Operators of Vector Calculus to Manifolds using Exterior Calculus}
\label{sec:ext_calc_ops}

The differential operators of vector calculus utilized in continuum mechanics formulations such as the $\mbox{grad}$, $\mbox{div}$, $\mbox{curl}$ can be extended to corresponding operators on general manifolds.  Differential operators on manifolds are notorious for having complicated notations when expressed in local coordinates~\cite{Abraham1988}.  We aim for a less coordinate-centric description of the methods and operators by utilizing approaches from exterior calculus.  For this purpose, we utilize the operators of exterior calculus given by the Hodge star $\star$, exterior derivative $\mb{d}$, and vector to co-vector isomorphisms $\flat, \sharp$ (definitions below). Operators extend to the context of general manifolds acting on scalar fields $f$ and vector fields $\mb{F}$ as
\begin{eqnarray}
\label{equ:gradDivCur}
\begin{array}{ll}
\mGrad_{\mathcal{M}}(f) = \lbrack \mb{d}f\rbrack^{\sharp}, &
\mDiv_{\mathcal{M}}(\mb{F}) =  -(-\star \mb{d}\star \mb{F}^\flat) = -\bs{\delta} \mb{F}^\flat, \\
\mCurl_{\mathcal{M}}(\mb{F}) =  -\star \mb{d}\left\lbrack \mb{F}^\flat \right\rbrack, &  \mCurl_{\mathcal{M}}(f) = \left\lbrack -\star \mb{d}f \right\rbrack^{\sharp} .
\end{array}
\end{eqnarray}
We define $\bs{\delta} = (-\star \mb{d}\star)$ which is referred to as the co-differential.  To define $\mb{d}$ the exterior derivative and $\star$ the Hodge star, we consider the tangent bundle $\mathcal{TM}$ of the manifold and its dual co-tangent bundle $\mathcal{TM}^*$.  The tangent bundle defines the spaces for scalar fields, vector fields, and more generally rank $m$ tensor fields over the manifold.  The co-tangent bundle is the space of duals to these fields.  The co-tangent bundle can be viewed as the space of differential forms of order $0$, $1$, and $m$. 

\qq{We denote vector fields and contravariant tensors using the notation $\mb{a} = a^{i_1\ldots i_k} \partial_{i_1}\cdots\partial_{i_k}$.} We use $\partial_{i_k}$ to denote the basis vector $\partial_{i_k} = \partial \bsy{\sigma}/\partial_{x^{i_k}}$ and tensor product these together to represent vectors and tensors for the choice of coordinates $\mb{x} = (x^1,x^2,\ldots,x^d)$.  We denote a differential $k$-form as 
$\bs{\alpha} = (1/k!)\alpha_{i_1,\ldots, i_k} \mb{d}\mb{x}^{i_1}\wedge \cdots \mb{d}\mb{x}^{i_k}$.  The $\wedge$ denotes the wedge-product of a tensor~\cite{Abraham1988}. We use the convention here with $1/k!$ to allow summations over all permutations of the index values for $i_1,\ldots,i_k$.  A more detailed discussion of tensor calculus on manifolds can be found in~\cite{Abraham1988}.  

We formulate the generalized operators in terms of the co-vectors (differential forms) $f^{\flat}$ and $\mb{F}^{\flat}$.  We use that in the case of a scalar field we have quantitatively at each point $f = f^{\flat}$~\cite{Abraham1988}.  
The isomorphisms $\flat,\sharp$ mapping between the vector and co-vector
spaces is given by 
\begin{eqnarray}
\label{eqn:op_iso_flat_sharp}
\mb{a}^{\flat}
& = & (1/k!) g_{i_1,\ell_1}\cdots g_{i_k,\ell_k} a^{\ell_1\ldots \ell_k} 
\mb{d}\mb{x}^{i_1}\wedge \cdots \mb{d}\mb{x}^{i_k} \\
\bs{\alpha}^{\sharp}
& = & (1/k!) g^{i_1,\ell_1}\cdots g^{i_k,\ell_k} \alpha_{\ell_1\ldots \ell_k}
\partial_{\mb{x}^{i_1}}\cdots\partial_{\mb{x}^{i_k}}.
\end{eqnarray}
The exterior derivative $\mb{d}$ of a differential $k$-form $\bs{\alpha}$ is defined in terms of the coordinates $\mb{x}$ as
\begin{eqnarray}
\label{eqn:op_ext_deriv}
\mb{d}\bs{\alpha} = \frac{1}{k!}\frac{\partial}{\partial x^j} \alpha_{i_1,\ldots, i_k} \mb{d}\mb{x}^{j} \wedge \mb{d}\mb{x}^{i_1}\wedge \cdots \mb{d}\mb{x}^{i_k}.
\end{eqnarray}
  The Hodge star $\star$ is defined in terms of the coordinates $\mb{x}$ as
\begin{eqnarray}
\label{eqn:op_hodge_star}
\star \bs{\alpha} = \frac{\sqrt{|g|}}{(n-k)!k!} \alpha^{i_1,\ldots,i_k} \epsilon_{i_1,\ldots,i_k,j_1,\ldots,j_{n-k}} \mb{d}x^{j_1}\wedge \cdots \wedge \mb{d}x^{j_{n-k}}.
\end{eqnarray}
Note the indices have been raised here for the $k$-form with $\alpha^{i_1,\ldots,i_k} = g^{i_1\ell_1}\cdots g^{i_k\ell_1}\alpha_{\ell_1,\ldots,\ell_k}$.  The $\epsilon_{\ell_1,\ldots,\ell_n}$ denotes the Levi-Civita tensor which gives the sign of the permutation of the indices $\ell_1,\ldots,\ell_n$ and is otherwise zero~\cite{Abraham1988}.  

\qq{This exterior calculus formulation allows us to provide a less coordinate centric description of the physics revealing in many cases more clearly the relationship of the continuum mechanics and role played by the geometry.  This also has the advantage in analytic calculations of making expressions more concise and allowing more readily for generalization of identities and techniques employed from vector calculation~\cite{AtzbergerGrossSurfExtCalc2017,AtzbergerSoftMatter2016}.}  As for practical numerical calculations, we  utilize  this  approach  along  with  symbolic  computation  to  generate offline the expressions needed for any choice of local coordinates on the manifold using equations~\ref{equ:gradDivCur}--~\ref{eqn:op_hodge_star}. This permits the efficient evaluation of these equations for any given choice of local coordinate using precompiled libraries for expressions. We give more details and show how this approach can be applied to the Laplace-Beltrami and Biharmonic operators in Appendix~\ref{appendix:monge_gauge}.

\subsubsection{GMLS Approximation of Differential Operators on Manifolds}
\label{sec:ext_calc_op_approx_manifolds}

\qq{Using these approaches we can perform GMLS estimates of the differential operators on the manifold.  We consider the approximation of target functionals which may depend nonlinearly on estimates of the geometry.  For example, the Laplace-Beltrami operator depends on the inverse metric tensor and can be expressed in local coordinates as}
\begin{equation}
\Delta_{\tiny LB} \phi = \frac{1}{\sqrt{|g|}} \partial_i \left( \sqrt{|g|} g^{ij} \partial_j \phi \right).
\end{equation}
\qq{We assume an estimate of $\mb{g}$ to be calculated at each point $\mb{x}_i$ following the approaches outlined in the previous sections.} We then approximate the action of the operator on scalar and vector fields through the following GMLS approach.  First, we find locally the best approximating reconstruction $\mb{P}$ of the scalar or vector field components on the manifold.  In the second, we apply the target functional for the differential operator to $\mb{P}$ using geometric quantities from our initial GMLS reconstruction of the manifold.  \qq{This can be expressed using the optimal coefficient vector $a_{\tilde{\mb{x}}}$ at $\tilde{\mb{x}}$ as}
\begin{equation}
\label{gmlsPairForm}
\tau^h_{\tilde{\mb{x}}}(\phi) = \tau_{\tilde{\mb{x}}}(\mathbf{P})^\intercal a_{\tilde{\mb{x}}}(u), \hspace{1cm}
a_{\tilde{\mb{x}}}(u) = \left(\Lambda(\mathbf{P})^\intercal \mathbf{W} \Lambda(\mathbf{P})\right)^{-1} \Lambda(\mathbf{P})^\intercal \mathbf{W} \Lambda(u).
\end{equation}
\qq{In the general setting, the sampling functionals $\lambda_j$ may depend non-linearly on the geometric terms.  In the current case using local point evaluations, the functionals $\lambda_j$ are linear.}

We remark that the two components $a_{\tilde{\mb{x}}}(u)$ and $\tau_{\tilde{\mb{x}}}(\mathbf{P})$ encode different types of information about the approximation.  The $\tau_{\tilde{\mb{x}}}(\mathbf{P})$ encodes the action of the target functional on the basis for the space $\mathbb{V}_h$.  The $a_{\tilde{\mb{x}}}(u)$ encodes the reconstruction of the function $u$ by the best approximating function $p^*$ in $\mathbb{V}_h$ according to the best match between the sampling functionals $\lambda_j$ acting on $u$ and $p^*$, see equation~\ref{eqn:opt_gmls_reconstruct}.  As a consequence, for each of the target operators $\tau$, the $a_{\tilde{\mb{x}}}(u)$ will not change since this term only depends on the function $u$.  As a result, we need only compute fresh for each operator the $\tau_{\tilde{\mb{x}}}(\mathbf{P})$ which represents how the differential operator on the manifold acts on the function space $\mathbb{V}_h$. \\ 

\noindent
\qq{As a summary, our GMLS approximation on the manifold involves the following steps}
\qq{
\begin{itemize}
\item Take $\Lambda = \{\lambda_j\}_{j=1}^N$ where $\lambda_j = \delta_{\mb{x}_j}$ are point evaluations with $\lambda_j\phi = \phi(\mb{x}_j)$ for the $\mb{x}_j$ in the neighborhood $j \in \mathcal{N}_i$ around the point $\mathbf{x}_i$.
\item For the target functionals $\tau$, treat the surface differential operators by utilizing for evaluation the parameterization and approximate metric tensor outlined in Section~\ref{sec:geo_approx_quant}.
\item For the reconstruction space $\mathbb{V}_h$, use the collection of $m_2^{th}$-order polynomials $p(x,y)$ over $\mathbb{R}^2$ where $m_2$ is an integer parameter for the maximum degree.
\item For  the weight function $\omega(\lambda_j,\tau_{\mb{x}_i}) = w(\|\mb{x}_j - \mb{x}_i \|)$, select a positive kernel $w(r)$ with support contained within an $\epsilon$-ball of $\mb{x}_i$.  We also use this to construct an $\epsilon$-graph on the point set.
\end{itemize}
}
\qq{The reconstruction space $\mathbb{V}_h$ consists of polynomials of order $m_2$ which need not be chosen to be the same order $m_1$ as in the geometric reconstructions in Section~\ref{sec:geo_reconstruct}, so in general $m_2 \neq m_1$.  However, in practice the operators on the manifold often involve differentiating geometric quantities which typically need $m_1 \geq m_2$ to achieve convergence. }

\qq{As an illustration of our GMLS approach, consider the Laplace-Beltrami operator.  The other differential operators for the manifold follow similarly, but with more complicated expressions which we evaluate symbolically, see Appendix~\ref{appendix:monge_gauge}.  For the Laplace-Beltrami operator we have
\begin{equation}
\Delta_{\tiny LB} \phi = \frac{1}{\sqrt{|g|}} \partial_i \left( \sqrt{|g|} g^{ij} \partial_j \phi \right),\hspace{1cm}
\label{eqn:LB_tau_rep}
\tau_{\tilde{\mb{x}}}(\mathbf{P};\mb{g}) = \frac{1}{\sqrt{|g|}} \partial_i \left( \sqrt{|g|} g^{ij} \partial_j \mathbf{P} \right).
\end{equation}
On the left the Laplace-Beltrami operator is expressed in coordinates and on the right we have the GMLS approximation.  The $\mb{P}$ represents the vector of basis functions of $\mathbb{V}_h$ and the differentials act component-wise.}

\begin{remark}  \qq{It is necessary to choose a reconstruction space $\mathbb{V}_h$ of sufficient richness that a differential operator on the manifold $\mathcal{L}_{\mb{g}}$ can be adequately represented.  For instance, a differential operator of order $k$ should have a polynomial space of order $m_2$ satisfying $m_2 \geq k$, as suggested by the bounds in equation~\ref{eqn:GMLS_deriv_bound}. From these bounds we also do expect 
that higher-order convergence rates are possible when using larger degrees. This can yield computational efficiencies in achieving a desired level of accuracy, especially when treating smooth fields and low-order differential operators.  However, it is important to note that larger choices of $m_2$ will necessitate that the neighborhoods defined by the weight function contain more points to ensure unisolvency and ultimately solvability of the GMLS problem.}
\end{remark}
We give additional details on our GMLS approach for specific operators in Appendix~\ref{appendix:monge_gauge}.

\section{Hydrodynamic Flows on Curved Surfaces}
\label{sec_governingEq}
\label{sec:hydro_formulation}

We formulate continuum mechanics equations for hydrodynamic flows on curved surfaces using approaches from the exterior calculus of differential geometry~\cite{Marsden1994,Abraham1988}.  This provides an abstraction that is helpful in generalizing many of the techniques of fluid mechanics to the manifold setting while avoiding many of the tedious coordinate-based calculations of tensor calculus.  The exterior calculus formulation also provides a coordinate-invariant set of equations helpful in providing insights into the roles played by the geometry in the hydrodynamics.  We provide a brief derivation of hydrodynamic equations here based on our prior work~\cite{AtzbergerSoftMatter2016,AtzbergerGrossSurfExtCalc2017,AtzbergerGrossHydroSurf2018}.   For additional discussion of the derivations for hydrodynamics on manifolds and related differential geometry, see~\cite{AtzbergerSoftMatter2016,AtzbergerGrossHydroSurf2018,SpivakDiffGeo1999,Marsden1994,Abraham1988}.


\subsection{Hydrodynamics in the Stokesian Regime}
\label{sec:hydro_formulation_stokes_regime}
We consider the hydrodynamics in the quasi-steady-state Stokes regime where the flow is determined by a balance between the fluid shear stresses and the body force.  The hydrodynamics in this regime can be expressed in covariant form as
\begin{eqnarray}
\label{equ_Stokes_geometric}
\label{equ_v_stokes}
\left\{
\begin{array}{llll}
\mu_m \left(-\bs{\delta} \mb{d} \mb{v}^{\flat} + 2 K \mb{v}^{\flat} \right) 
- \gamma \mb{v}^{\flat} - \mb{d}p 
& = & - \mb{b}^{\flat} \\
-\bs{\delta} \mb{v}^{\flat} & = & 0.
\end{array}
\right.
\end{eqnarray}
The $\mb{v}^{\flat}$ is the surface fluid velocity, $p$ the surface presssure enforcing incompressibility, and $\mb{b}^{\flat}$ the surface force density driving the flow.  The $\mu_m \left(-\bs{\delta} \mb{d} \mb{v}^{\flat} + 2 K \mb{v}^{\flat} \right)$ corresponds to the divergence of the internal shear stress of the surface fluid, and $-\bs{\delta} \mb{v}^{\flat} =  0$ expresses the incompressibility constraint.  The $\mu_m$ gives the surface fluid viscosity.  It is worth pointing out that the surface shear stress has a dependence not only on the usual gradients in the velocity field but also the Gaussian Curvature $K$ of the surface.  This can lead to interesting flow phenomena on curved surfaces and significant differences with respect to flat surfaces, as discussed in~\cite{AtzbergerSoftMatter2016,AtzbergerGrossHydroSurf2018,ArroyoRelaxationDynamics2009,LevineDinsmoreHydroEffectTopology2008}.  

We remark that the $-\gamma\mb{v}^{\flat}$ serves as our model for the coupling between the surface flow and bulk three-dimensional surrounding fluid.  More sophisticated models also can be formulated, but for general geometries this requires development of a separate solver for the bulk three-dimensional surrounding fluid which we shall consider in future work.  It is important in physical models to have some form of dissipative traction stress with the surrounding bulk fluid since this provides a crucial dissipative mechanism that suppresses the otherwise well-known Stokes paradox that arises in purely two-dimensional fluid equations~\cite{Acheson1990,batchelor2000,AtzbergerGrossHydroSurf2018,Seki2007,Saffman1976}.  Additional discussions of equation~\ref{equ_Stokes_geometric} and its derivation can be found in~\cite{AtzbergerGrossHydroSurf2018,AtzbergerSoftMatter2016}.

\subsection{Vector Potential Formulation for Incompressible Flows and Hodge Decomposition}
\label{sec:hydro_formulation_vec_potential}
We generalize approaches from fluid mechanics to the context of manifolds to handle the incompressibility constraint in equation~\ref{equ_Stokes_geometric}.  We reformulate equation~\ref{equ_Stokes_geometric} using the Hodge decomposition and a vector potential $\phi$ that ensures the generated velocity fields are incompressible.  By utilizing this gauge to describe the physics we can avoid the challenges in numerical methods associated with having to enforce explicitly the incompressibility constraint.  We use a surface Hodge decomposition of the fluid velocity field that can be expressed using the exterior calculus as
\begin{eqnarray}
\label{equ_hodge_decomp}
\mb{v}^{\flat} = \mb{d}\psi + \bs{\delta}\phi + \mb{h}.
\end{eqnarray}
The $\psi$ is a $0$-form, $\phi$ is a $2$-form, and $\mb{h}$ is a harmonic $1$-form on the surface with respect to the Hodge Laplacian $\Delta_H \mb{h} = \left(\bs{\delta}\mb{d} + \mb{d}\bs{\delta}\right) \mb{h} = 0$.  The first term $\mb{d}\psi$ captures the curl-free component of the velocity field, the second term $\bs{\delta}\phi$ the divergence-free component of the velocity field, and the third term an additional harmonic part that arises from the topology of the manifold.  In the Euclidean setting only the first two terms typically play a role since the harmonic term in this case is often a trivial constant and with decay conditions at infinity the constant is zero.  

In the non-Euclidean setting there can be many non-trivial harmonic $1$-forms.  The number is determined by the dimensionality of the null-space of the Hodge Laplacian which depends on the topology of the manifold~\cite{JostBookDiffGeo1991}.  As a consequence, we have for different topologies that the richness of the harmonic differential forms $\mb{h}$ appearing in equation~\ref{equ_hodge_decomp} will vary.  Fortunately, in the case of spherical topology the surface admits only the trivial harmonic $1$-forms $\mb{h} = \bs{0}$ making this manifold relatively easy to deal with in our physical descriptions.  As we shall discuss, for more general topologies our incompressibility gauge descriptions will require solving additional coupled equations in order to resolve the non-trivial harmonic contributions.  We shall focus here primarily on the case of manifolds having spherical topology and pursue in future work development of these additional numerical solvers needed for the harmonic component.

We consider incompressible velocity fields $\mb{v}^{\flat}$ on manifolds having spherical topology.  When applying the co-differential $\bsy{\delta}$ to equation~\ref{equ_hodge_decomp} and utilizing the incompressibility constraint in equation~\ref{equ_Stokes_geometric}, we have $\bs{\delta} \mb{v}^{\flat} = \bs{\delta}\mb{d} \psi = \Delta_H \psi = 0$.  For spherical topology this requires $\psi = C$ and $\mb{d}\psi = 0$.  As a consequence, we can express the incompressible hydrodynamic velocity fields as
\begin{eqnarray}
\label{equ_hodge_decomp2}
\mb{v}^{\flat} = \bs{\delta}\phi.
\end{eqnarray}
From the co-differential operator $\bsy{\delta}$ defined in Section~\ref{sec:ext_calc_ops}, we see that $\phi$ is a $2$-form on the two-dimensional surface.  In practice, we find it more convenient to express $\mb{v}^{\flat}$ in terms of an operation on a $0$-form (scalar field) which can be done using the Hodge star to obtain $\Phi = \star \phi$.  Using the identity of the Hodge star that $\star\star = (\star)^2 = -1$ for 2-manifolds.  This gives $\phi = -\star \Phi$.  This allows us to express incompressible hydrodynamic flow fields as 
\begin{eqnarray}
\label{eqn_v_curl_rep}
\mb{v}^{\flat} = -\star\mb{d} \Phi.
\end{eqnarray}
\qq{We refer to $\Phi$ as the \textit{vector potential} since it serves as a potential to generate vector fields $\mb{v}$.  This can be interpreted as a generalized curl operation as in equation~\ref{equ:gradDivCur} applied to a scalar field which intrinsically generates divergence-free vector fields $\mb{v}$ .  This approach generalizes the vorticity-stream formulation of fluid mechanics~\cite{Acheson1990a} to the manifold setting.   We use this to reformulate the hydrodynamic equations in terms of unconstrained equations in terms of $\Phi$.}

\subsubsection{Biharmonic Formulation of the Hydrodynamics}
\label{sec:hydro_formulation_biharmonic}
\qq{We reformulate the hydrodynamics equations~\ref{equ_v_stokes} in terms of an unconstrained equation for the vector potential $\Phi$.}  We substitute equation~\ref{eqn_v_curl_rep} into equation~\ref{equ_v_stokes} and apply the generalized curl operator $\mbox{curl}_\mathcal{M} = -\star\mb{d}$ to both sides.  This gives the biharmonic hydrodynamic equations on the surface
\begin{eqnarray}
\label{equ_unsplit_stokes}
-\mu_m\Delta_H^2 \Phi - \gamma \Delta_H \Phi - 2\mu_m (-\star\mb{d}(K(-\star\mb{d})))\Phi &=& -\star\mb{d} \mb{b}^{\flat}.
\end{eqnarray}
The $\mu_m$ is the surface shear viscosity, $\gamma$ the drag with the surrounding bulk fluid, and $K$ the Gaussian curvature of the manifolds.  The $\mb{b}^{\flat}$ is the covariant form for the body force acting on the fluid.  We see the pressure term no longer plays a role relative to equation~\ref{equ_v_stokes}.

The Hodge Laplacian now acts on $0$-forms as $\Delta_H \Phi = \bs{\delta}\mb{d}\Phi$ and is related the surface Laplace-Beltrami operator by $\Delta_H \Phi = -\Delta_{LB}\Phi$.  This provides for numerical methods a particularly convenient form for the fluid equations since it only involves solving for a scalar field $\Phi$ on the surface.  However, this does have the drawback that for handling the incompressibility constraint this way we now need to solve a biharmonic equation on the surface.  We shall refer in our numerical methods to this approach to the hydrodynamics as the \textit{biharmonic formulation}.

We remark that our approach can be related to classical methods in fluid mechanics by viewing our operator $-\star\mb{d}$ as a type of curl operator that is now generalized to the manifold setting.  The $\Phi$ serves the role of a vector potential for the flow~\cite{Acheson1990,batchelor2000,Lamb1895}.  The velocity field of the hydrodynamic flows $\mb{v}$ is recovered from the vector potential $\Phi$ as $\mb{v}^{\flat} = -\star \mb{d} \Phi$.  We obtain the velocity field $\mb{v} = \mb{v}^{\sharp} = \left(-\star \mb{d} \Phi\right)^{\sharp}$ using equation~\ref{equ_gen_curl_0form} and the isomorphisms $\sharp$ between co-vectors and vectors discussed in Section~\ref{sec:ext_calc_ops}.  Additional discussion of this formulation of the hydrodynamics can be found in~\cite{AtzbergerSoftMatter2016,AtzbergerGrossHydroSurf2018}.  

\subsubsection{Split Formulation of the Hydrodynamics}
\label{sec:hydro_formulation_split}
While the equation~\ref{equ_unsplit_stokes} is expressed in terms of biharmonic operators, for numerical purposes we can reformulate the problem by splitting it into two sub-problems each of which only involve the Hodge Laplacian.  This is helpful since for our numerical methods this would require us to only need to resolve second order operators with our GMLS approximations.  This has the practical benefit of greatly reducing the size of the GMLS stencil sizes ($\epsilon$-neighborhoods) required for unisolvency for the operator as discussed in Section~\ref{sec:gmls}. 

We reformulate the hydrodynamic equations by defining $\Psi = \Delta_H \Phi$, which allows us to split the action of the fourth-order biharmonic operator into two equations involving only second- order Hodge Laplacian operators as
\begin{eqnarray}
\label{equ_split_stokes}
-\mu_m\Delta_H \Psi - \gamma \Psi - 2\mu_m (-\star\mb{d}(K(-\star\mb{d})))\Phi &=& -\star\mb{d} \mb{b}^{\flat}. \\
\Delta_H \Phi - \Psi &=& 0.
\end{eqnarray}
\qq{As we shall discuss, the lower order of the differentiation has a number of benefits even though we incur the extra issue of dealing with a system of equations.}  This reformulation results in less sensitivity to errors in the underlying approximations in the GMLS reconstructions of the geometry and surface fields.  This reformulation also requires much less computational effort and memory when assembling the stiffness matrices since the lower order permits use of smaller $\epsilon$-neighborhoods to achieve unisolvency as discussed in Section~\ref{sec:gmls}.  We refer to this reformulation of the hydrodynamic equations as the $\textit{split formulation}$.

For a further discussion of these surface hydrodynamics equations, related derivations, and physical phenomena see~\cite{AtzbergerSoftMatter2016,AtzbergerGrossHydroSurf2018}.

\section{Computational Methods and Numerical Solvers}
\label{sec:numerical_solver}

\begin{figure}[H]
\centerline{\includegraphics[width=0.95\columnwidth]{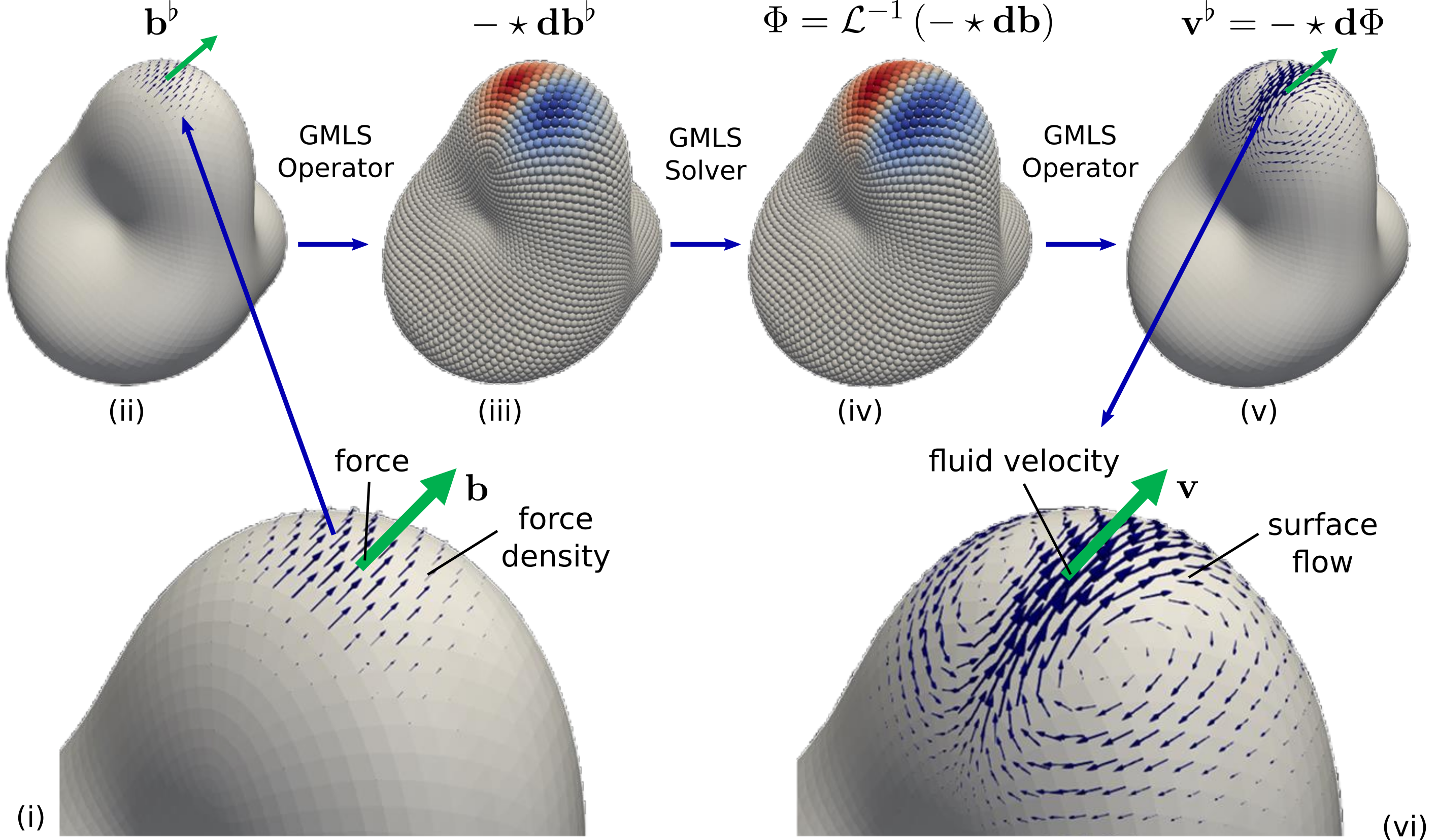}}
\caption{\qq{Steps in the Numerical Methods for the Surface Hydrodynamic Flows.  For a given body force density or stresses $\mb{b}$ acting on the surface fluid we convert the fields to covariant form $\mb{b}^{\flat}$, shown in (i),(ii).  To handle incompressibile flows, we convert all fields to a divergence-free gauge using the generalized surface curl $-\star \mb{d} \mb{b}^{\flat}$, shown in (iii).  We solve for the vector potential $\Phi$ of the surface hydrodynamic flow using equations~\ref{equ_unsplit_stokes} or~\ref{equ_split_stokes} and our GMLS collocation methods for the differential operators, shown in (iv).  We construct the covariant form of the velocity field of the hydrodynamic flow response using the generalized surface curl $\mb{v}^{\flat}$, shown in (v).  We obtain our final results by converting the covariant form $\mb{v}^{\flat}$ to the velocity field by $\mb{v} = \left(\mb{v}^{\flat}\right)^{\sharp}$.  This yields the surface hydrodynamic flow shown in (vi).}}
\label{fig:solver_schematic}	
\end{figure}

\qq{We develop numerical methods to solve equations~\ref{equ_unsplit_stokes} or~\ref{equ_split_stokes} 
for the velocity field of hydrodynamic flows on surfaces using the GMLS approximations of Section~\ref{sec_gmls} and~\ref{sec:geo_reconstruct}.  We briefly discuss the overall steps used in our numerical methods.  We formulate the hydrodynamics using a vector-potential formulation to obtain a gauge that intrinsically enforces the incompressibility constraints of the flow appearing in equation~\ref{equ_v_stokes}.}  For steady-state hydrodynamic flows, we derived conditions for the vector potential of the flow resulting in equation~\ref{equ_unsplit_stokes}.  We summarize the steps used in our solution approach in Figure~\ref{fig:solver_schematic}.

\qq{To determine numerically the hydrodynamic flow in response to a body force density $\mb{b}$ acting on the surface fluid, we first convert force fields into co-variant form $\mb{b}^{\flat}$.}  We next use our exterior calculus formulation of the generalized curl to obtain the corresponding vector-potential for the body force $\Psi = \mathcal{C}_1 \mb{b}^{\flat}$ where $\mathcal{C}_1 = -\star \mb{d}$ acts on $1$-forms.  We numerically compute $\Psi = \tilde{C}_1 \mb{b}^{\flat}$ where $\tilde{C}_1$ is our GMLS approximation of the curl operator $\mathcal{C}_1 $ discussed in Section~\ref{sec:ext_calc_op_approx_manifolds} and Appendix~\ref{appendix:monge_gauge_coord_ext_calc_exp}.

We can now utilize equation~\ref{equ_unsplit_stokes} to specify the differential equation for the steady-state velocity response.  We use GMLS to assemble in strong form a stiffness matrix $A$ using a collocation approach.  The full differential operator that appears on the left-hand-side is computed at each base point $\tilde{\mb{x}}$ of the point set of the manifold.  This results in the system of equations linear in $\tilde{\mb{\Phi}}$
\begin{eqnarray}
\label{eqn:linear_sys_hydro}
A \tilde{\mb{\Phi}} = \tilde{C}_0 \mb{b}.
\end{eqnarray}
We solve the large linear system using GMRES with algebraic multigrid (AMG) preconditioning.

The velocity field is given from the vector potential $\Phi$ by the generalized surface curl operator $\mb{v}^{\flat} = \mathcal{C}_0 \Phi$, where $\mathcal{C}_0 = -\star \mb{d}$ acts on $0$-forms.  From the solution $\tilde{\mb{\Phi}}$ of equation~\ref{eqn:linear_sys_hydro}, we construct numerically the co-variant velocity field of the flow using $\tilde{\mb{v}}^{\flat} = \tilde{C}_0 \tilde{\mb{\Phi}}$.  The $\tilde{C}_0$ is our GMLS approximation of the generalized curl operator $\mathcal{C}_0$ discussed in Section~\ref{sec:ext_calc_ops}.  Finally, using the metric tensor obtained from the GMLS reconstruction, we obtain the surface velocity field $\tilde{\mb{v}}$ by converting the covariant field $\mb{v}^{\flat}$ into the contravariant field by $\tilde{\mb{v}} = \left(\tilde{\mb{v}}^{\flat}   \right)^{\sharp}$.  For more details on this approach and operations see Appendix~\ref{appendix:monge_gauge_coord_ext_calc_exp}.  We use this approach to numerically compute incompressible hydrodynamic flows in response to applied driving forces or stresses acting on the surface fluid.  \qq{We remark that our approach can also be combined with other computational methods and solvers to compute coupling to bulk three dimensional hydrodynamics or more generally for resolving other physical systems and interactions that occur at interfaces.}

\qq{The Compadre toolkit~\cite{Compadre2019} was used to solve the GMLS problems.  The toolkit provides domain decomposition and distributed vector representation of fields as well as global matrix assembly. The linear equations were solved through iterative block solvers (Belos~\cite{Belos}), block preconditioners (Teko) and AMG preconditioning (MueLu~\cite{MueLu2014a,MueLu2014b}), within the Trilinos software framework~\cite{trilinos2005}.}

\section{Results}
\label{sec_results}
\label{sec:results}

\subsection{Convergence Results for Operators on Manifolds based on GMLS Geometric Reconstructions}
\label{sec:conv_ops}

We investigate the convergence of the GMLS approximation of the operators required to solve the hydrodynamic equations in Section~\ref{sec:hydro_formulation}.  It is important to note that   
our target functionals have a non-linear dependence on the geometry resulting in contributions from two different GMLS approximations.  First, the GMLS reconstruction of the geometry of the manifold and associated associated geometric quantities.  Second, the GMLS approximation of differential operators acting on the surface scalar and vector fields.

Our solvers for the surface hydrodynamic flows of equation~\ref{equ_unsplit_stokes} and~\ref{equ_split_stokes} use the following operators: \textit{Laplace-Beltrami} $\mathcal{L}_{LB} = -\Delta_H= -\mb{d}\bsy{\delta}$, \textit{Biharmonic} $\mathcal{L}_{BH} = \mathcal{L}_{LB}^2 = \Delta_H^2$,\textit{Curvature} $\mathcal{L}_{K} = curl_{\mathcal{M}} \left(K \cdot curl_{\mathcal{M}} \right) = -\star \mb{d}\left(K \cdot (-\star\mb{d}\right)$, \textit{Surface-Curl-$0$-Forms $\mathcal{C}_0$},$\mathcal{L}_{C0} = \mbox{curl}_\mathcal{M}\Phi = -\star\mb{d}_0 $, \textit{Surface-Curl-$1$-Forms $\mathcal{C}_1$},$\mathcal{L}_{C1} = \mbox{curl}_\mathcal{M} \mb{v} = -\star\mb{d}_1$.  In the split formulation of equation~\ref{equ_split_stokes}, this simplifies without the need for $\mathcal{L}_{BH}$.

\qq{To approximate each of the operators using GMLS, we use point samples within a distance $\epsilon$ from the target point $x_i$.  For the local optimization problems to determine the polynomial of degree $m$, a minimum of $n_p = \binom{m+2}{2}$ points are needed per neighborhood.  To find such an $\epsilon$, we sweep through the points and either increase or decrease the initial guess $\epsilon_0$ by a factor $\beta_*$.  This is done so that the $\epsilon = \beta_*^k\epsilon_0$ can be used globally to generate neighborhoods with at least $\alpha_*n_p$ points in each of them, where $\alpha_* \geq 1$ is a tune-able parameter for calculations.  This is used to determine the polynomial representation $p^*$.  In practice, the factor $\alpha_*$ is used to ensure significantly more points are in neighborhoods than are strictly needed which helps to provide additional robustness in calculations.  We use throughout $\alpha_* = 2.8$ and $\beta_* = 2$.
}

To study the accuracy of our GMLS approximation of these operators, we investigate the case of the test scalar field $\Phi(\mb{X}) = \Phi(x,y,z) = z(x^4 + y^4 - 6x^2y^2)$ and test vector field $\mb{v}^{\flat} = \mathcal{C}_0\Phi = -\star \mb{d}_0 \Phi$.  We have chosen $\Phi(x,y,z)$ to be a smooth continuation of a spherical harmonic mode to the full space $\mathbb{R}^3$.   Since our manifolds $\mathcal{M}$ are smooth, we obtain a smooth surface scalar field by evaluation of $\Phi(\mb{X})$ at the surface (i.e. using the inclusion map $\iota: \mathbb{R}^3 \xhookrightarrow{} \mathcal{M}$ with $\Phi(\mb{x}) = \iota_{\mb{x}} \Phi(\cdot)$).  This provides a way to define surface scalar fields and vector fields without the need for local coordinate charts on the manifold.  

We investigate the accuracy of the GMLS approximation of these operators.  We study the $\ell^2$-errors
\begin{eqnarray}
\epsilon_{op0} = \left\|\tilde{L}_{\mb{g}} \Phi - \mathcal{L}_{\mb{g}} \Phi \right\|_2, \hspace{1cm}
\epsilon_{op1} = \left\|\tilde{L}_{\mb{g}} \mb{v} - \mathcal{L}_{\mb{g}} \mb{v} \right\|_2.
\end{eqnarray}
The $\ell^2$-norm is computed by averaging the error over all $n$ sample points of the manifold $\|u - v\|_2^2 = \frac{1}{n} \sum_i \left(u(\mb{x}_i) - v(\mb{x}_i)\right)^2$.
The $\tilde{L}_{\mb{g}}$ denotes the numerical GMLS approximation of the operator $\mathcal{L}_{\mb{g}}$.  \qq{In practice, for comparison with the GMLS results in the convergence studies, we evaluate to high precision the action of the operators $\mathcal{L}_{\mb{g}}$ by using symbolic calculations using SymPy~\cite{Sympy2017}.}  

Using this approach, we investigate the accuracy of the GMLS approximation of the operators for each of the manifolds in Tables~\ref{table:conv_L_LB}--~\ref{table:conv_L_K}.  We estimate approximate convergence rates by fitting using in a log-log scale the error between the reported $h$ value and the previous $h$ value.  \qq{While we do not have theory given that the operators have a non-linear dependence on the manifold geometry, for operators of order $k$ and GMLS approximation of order $m$ we do have the suggestive predictions from equation~\ref{eqn:GMLS_deriv_bound} that the convergence be order $m + 1 - k$.}  Since our GMLS methods involve approximations both of the geometry and the surface fields, for purposes of the comparisons we take $k = \max(k_1,k_2)$ and $m = \max(m_1, m_2)$.  The
$k_1$ denotes the order of the differentiation involved in obtaining the quantities associated with the geometry and $k_2$ with the order of differentiation of the surface fields.  The $m_1, m_2$ are the polynomial orders used for the approximations for the manifold geometry and surface fields, as discussed in Section~\ref{sec:gmls}.

\begin{table}[H]
\label{LBConv}
\begin{center}
\begin{tabular}{ccccccccccc}
\textbf{}  & \multicolumn{2}{c}{\textbf{Manifold A}}      & \multicolumn{2}{c}{\textbf{Manifold B}} & \multicolumn{2}{c}{\textbf{Manifold C}} & & & \multicolumn{2}{c}{\textbf{Manifold D}}      \\\cline{2-7}\cline{10-11}
\textbf{h} & \textbf{$\ell_2$-error} & \textbf{Rate} & \textbf{$\ell_2$-error} & \textbf{Rate} & \textbf{$\ell_2$-error} & \textbf{Rate} & & \textbf{h} & \textbf{$\ell_2$-error} & \textbf{Rate} \\
0.1 & 4.2208e-04 & - & 2.2372e-02 & -  & 1.3580e-01  & -&  & .08 &  4.7880e-02 & -  \\
0.05 & 7.503e-06 & 5.74 & 1.2943e-03 & 4.11 & 4.8597e-03 & 4.80 & & .04 & 5.5252e-04 & 6.54 \\
0.025 & 1.8182e-07 & 5.34  & 5.8300e-05 & 4.46 & 1.2928e-04 & 5.24 & & .02 &  1.3877e-05 & 5.36\\
0.0125 & 4.8909e-09 & 5.21 & 1.7364e-06 & 5.06 & 3.7508e-06 &  5.11 & & .01 & 3.7568e-07 & 5.17 \\
\end{tabular}
\end{center}
\caption{Convergence of GMLS Approximation of the Laplace-Beltrami Operator $\mathcal{L}_{LB}$.  We use GMLS with ($k = 2, m = 6)$ and find the methods have $\sim 5^{th}$-order asymptotic convergence.  The target sampling distance $h$ is discussed in Appendix~\ref{appendix:sampling_res_manifolds}.}
\label{table:conv_L_LB}
\end{table}

\begin{table}[H]
\label{BiLBConv}
\begin{center}
\begin{tabular}{ccccccccccc}
\textbf{}  & \multicolumn{2}{c}{\textbf{Manifold A}}      & \multicolumn{2}{c}{\textbf{Manifold B}} & \multicolumn{2}{c}{\textbf{Manifold C}} & & & \multicolumn{2}{c}{\textbf{Manifold D}}      \\\cline{2-7}\cline{10-11}
\textbf{h} & \textbf{$\ell_2$-error} & \textbf{Rate} & \textbf{$\ell_2$-error} & \textbf{Rate} & \textbf{$\ell_2$-error} & \textbf{Rate} & & \textbf{h} & \textbf{$\ell_2$-error} & \textbf{Rate} \\
0.1 & 1.7177e-01 & - & 1.1102e+01 & -  & 6.9226e+01  & -&  & .08 &  4.0566e+01 & -  \\
0.05 & 1.0768e-02 & 3.94 & 2.1455e+00 & 2.37 & 9.6017e+00 & 2.85 & & .04 & 1.3004e+01 & 5.04 \\
0.025 & 9.3281e-04 & 3.51 & 3.4556e-01 & 2.63 & 7.8738e-01 & 3.61 & & .02 &  1.0736e-01 & 3.63 \\
0.0125 & 9.3585e-05 & 3.31 & 3.5904e-02 & 3.26 & 7.7925e-02 &  3.34 & & .01 & 1.0722e-02 & 3.30 \\
\end{tabular}
\end{center}
\caption{Convergence of GMLS Approximation of the Biharmonic Laplace-Beltrami Operator $\mathcal{L}_{BH} = \mathcal{L}_{LB}^2$.  We use GMLS with ($k = 4, m = 6)$ and find the methods have $\sim 3^{rd}$-order asymptotic convergence.  }
\label{table:conv_L_BH}
\end{table}

\begin{table}[H]
\label{curlKcurlConv}
\begin{center}
\begin{tabular}{ccccccccccc}
\textbf{}  & \multicolumn{2}{c}{\textbf{Manifold A}}      & \multicolumn{2}{c}{\textbf{Manifold B}} & \multicolumn{2}{c}{\textbf{Manifold C}} & & & \multicolumn{2}{c}{\textbf{Manifold D}}      \\\cline{2-7}\cline{10-11}
\textbf{h} & \textbf{$\ell_2$-error} & \textbf{Rate} & \textbf{$\ell_2$-error} & \textbf{Rate} & \textbf{$\ell_2$-error} & \textbf{Rate} & & \textbf{h} & \textbf{$\ell_2$-error} & \textbf{Rate} \\
0.1 & 3.7004e-03 & - & 1.0621e+01 & -  & 6.1440e+01  & -&  & .08 &  6.5445e-01 & -  \\
0.05 & 1.9863e-04 & 4.16 & 1.7987e-01 & 2.56 & 3.9161e-01 & 3.97 & & .04 & 1.6209e-02 & 5.42 \\
0.025 & 1.1937e-05 & 4.03 & 1.9796e-02 & 3.18 & 2.9043e-02 & 3.76 & & .02 &  8.4581e-04 & 4.30 \\
0.0125 & 7.3369e-07 & 4.01 & 1.6147e-03 & 3.61 & 2.0897e-03 &  3.80 & & .01 & 5.6742e-05 & 3.87 \\
\end{tabular}
\end{center}
\caption{Convergence of GMLS Approximation of the Curl-K-Curl Operator $\mathcal{L}_{K}$.  We use GMLS with ($k_1 = 3, k_2 = 2, m = 6)$ and find the methods have $\sim 4^{rd}$-order asymptotic convergence.  }
\label{table:conv_L_K}
\end{table}

We find to a good approximation our GMLS methods exhibit convergence rates in agreement with the suggestive prediction $m + 1 - k$.  For the Laplace-Beltrami operator $\mathcal{L}_{LB}$ with $(k = 2, m = 6)$, we find $\sim 5^{th}$-order convergence rate, see Table~\ref{table:conv_L_LB}.  For the Biharmonic operator $\mathcal{L}_{BH}$ with $(k = 4, m = 6)$, we find $3^{rd}$-order convergence rate, see Table~\ref{table:conv_L_BH}.  In the case of the Curvature Operator $\mathcal{L}_K$ we have $(k_1 = 3, k_2 = 2, m = 6)$.  The $k_1 = 3$ arises since the operator involves estimation not only of the surface Gaussian Curvature $K$ but also its first derivatives.  For $\mathcal{L}_K$, we find $\sim$ 4th-order convergence rate, see Table~\ref{table:conv_L_K}.  We also report convergence rates for the curl operators $\mathcal{L}_{C0}$ and $\mathcal{L}_{C1}$ in Appendix~\ref{appendix:conv_curl}. \qq{We give further details on the sampling resolution of the manifolds in Appendix~\ref{appendix:sampling_res_manifolds}.}  \qq{We perform further convergence studies to investigate the robustness of the methods and how the accuracy depends on the quality of the point sampling of the manifold geometry in Appendix~\ref{appendix:stability_sampling}.  Again, we emphasize while there is currently no rigorous convergence theory given the non-linear dependence on geometry in our GMLS approximations, we do find in each case agreement with the suggestive predictive rates $m + 1 - k$ similar to equation~\ref{eqn:GMLS_deriv_bound}.  } 

\subsection{Convergence Results for Hydrodynamic Flows}
\label{sec:stokes_conv}

\begin{figure}[H]
\centerline{\includegraphics[width=0.8\columnwidth]{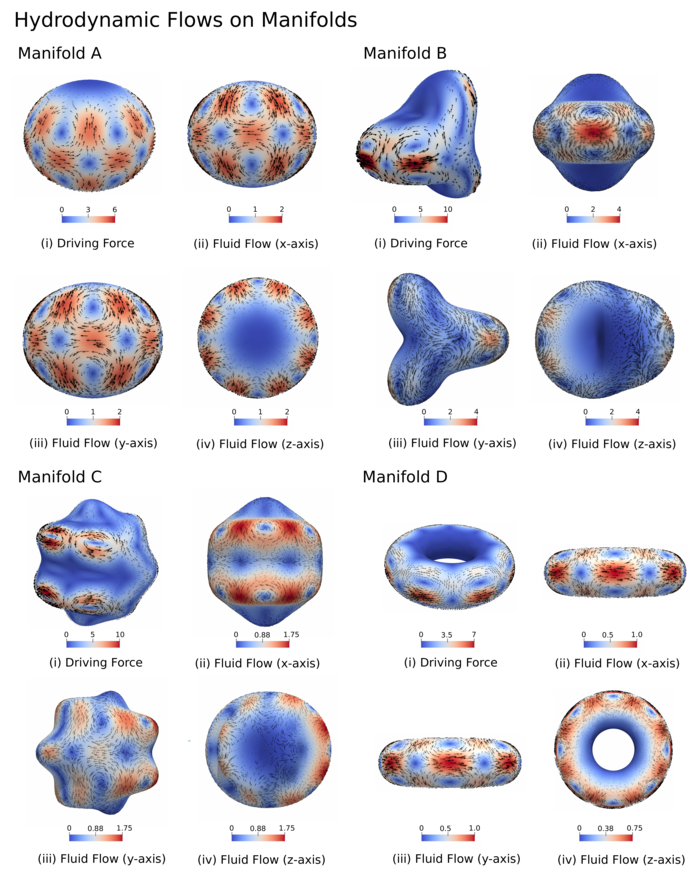}}
\caption{\qq{Surface Hydrodynamic Flows on Manifolds $A$--$D$.  We use our GMLS solver to compute numerically the surface hydrodynamic flow responses on each of the manifolds, as discussed in Section~\ref{sec:hydro_formulation_vec_potential}.  Manifold $A-D$ solutions were computed at a resolution with number of sample points
$n_A = 38,486, n_B = 147,634, n_C = 127,346,$ and $n_D = 118,942$.}}
\label{fig:hydro_surf_flows_A_B}
\label{fig:hydro_surf_flows_C_D}	
\label{fig:hydro_surf_flows_all}	
\end{figure}

We investigate the convergence of our GMLS methods for the surface hydrodynamic equations formulated in Section~\ref{sec:hydro_formulation}.  We study convergence of our solvers for hydrodynamic flows by developing manufactured solutions using high precision symbolic calculations of the incompressible flow field $\mb{v}^{\flat} = -\star \mb{d} \Phi = \mathcal{C}_0 \Phi$ with the specific choice of $\Phi$ given in Section~\ref{sec:conv_ops}.  

We calculate symbolically the expressions of the forcing term $\mb{b}$ using  equation~\ref{equ_Stokes_geometric} where $\mu_m(-\bsy{\delta} \mb{d} + 2K) \mb{v}^\flat - \gamma \mb{v}^\flat - \mb{d}p = -\mb{b}^\flat$.  We manufacture the data $\mb{b}$ needed on the RHS of equation~\ref{equ_Stokes_geometric} using
\begin{equation}
\label{eqn:manufacture_b}
    \mb{b}^{\flat} = \mu_m \bsy{\delta}\mb{d} \mb{v}^{\flat} + (\gamma- 2\mu_m K) \mb{v}^{\flat}.
\end{equation}
Since generating both the velocity field $\mb{v}$ and force density $\mb{b}$ this way will already be incompressible, we have used that we can set $p = 0$ when manufacturing our data.  In practice, we evaluate equation~\ref{eqn:manufacture_b} to high precision using the symbolic package SymPy~\cite{Sympy2017}.

We investigate the convergence of the GMLS solvers using the $\ell^2$-error
\begin{eqnarray}
\epsilon_{hydro} & = & \left \|\tilde{\mb{v}} - \mb{v} \right \|_2/ \| \mb{v} \|_2 \\
\tilde{\mb{v}} & = & \mathcal{C}_0 (\mathcal{S}^{-1} (\mathcal{C}_0 \mb{b})).
\end{eqnarray}
The $\mb{v}$ denotes the exact solution, $\mathcal{C}_0$ approximates numerically $-\star \mb{d}_0$, $\mathcal{C}_1$ approximates numerically $-\star \mb{d}_1$, and $\mathcal{S}^{-1}$ denotes the numerical solution operator corresponding to use of our GMLS solver.  We use the hydrodynamics equations both formulated using the biharmonic form in equation~\ref{equ_unsplit_stokes} or in the split form in equation~\ref{equ_split_stokes}.

For each of the manifolds $A-D$, we computed manufactured solutions with the parameters $\mu_m = 0.1$, $\gamma = 0.1$ in equation~\ref{eqn:manufacture_b}.  We used the surface force density $\mb{b}$ to numerically compute surface hydrodynamic flow responses $\tilde{\mb{v}}$ using our GMLS solvers discussed in Section~\ref{sec:numerical_solver}.  We show the hydrodynamic surface flows in Figure~\ref{fig:hydro_surf_flows_all}.  \qq{We show our convergence results for both the case of the biharmonic formulation and split formulation in the Tables~\ref{table:hydro_conv_bh_and_split_A}--~\ref{table:hydro_conv_bh_and_split_D}.
}

\begin{table}[H]
\begin{center}
Biharmonic Formulation: Manifold A\\
\begin{tabular}{ccccccc}
\textbf{}  & \multicolumn{2}{c}{\textbf{m = 4}}      & \multicolumn{2}{c}{\textbf{m = 6}}      & \multicolumn{2}{c}{\textbf{m = 8}}      \\\cline{2-7}
\textbf{h} & \textbf{$\ell_2$-error} & \textbf{Rate} & \textbf{$\ell_2$-error} & \textbf{Rate} & \textbf{$\ell_2$-error} & \textbf{Rate} \\
0.1    & 1.6072e-01 & -     & 1.1597e-03 & -    & 1.0648e-03 & -    \\
0.05   & 1.8027e-02 & 3.11 & 8.4190e-05 & 3.73 & 1.8627e-06 & 9.04 \\
0.025  & 4.9155e-03 & 1.86 & 1.1655e-05 & 2.84 & 4.4796e-08 & 5.35 \\
0.0125 & 2.0873e-03 & 1.23  & 7.1161e-07 & 4.02 & 1.9263e-07 & -2.10 \\
\end{tabular}
\\
\vspace{0.2cm}
Split Formulation: Manifold A\\
\begin{tabular}{ccccccc}
\textbf{}  & \multicolumn{2}{c}{\textbf{m = 4}}      & \multicolumn{2}{c}{\textbf{m = 6}}      & \multicolumn{2}{c}{\textbf{m = 8}}      \\\cline{2-7}
\textbf{h} & \textbf{$\ell_2$-error} & \textbf{Rate} & \textbf{$\ell_2$-error} & \textbf{Rate} & \textbf{$\ell_2$-error} & \textbf{Rate} \\
0.1    & 1.5578e-02 & -     & 2.6826e-04 & -    & 1.0756e-04 & -    \\
0.05   & 7.0783e-04 & 4.40     & 1.2065e-05 & 4.41 & 3.7309e-07 & 8.06 \\
0.025  & 1.2151e-05 & 5.83  & 4.4532e-07 & 4.74 & 3.0556e-09 & 6.90 \\
0.0125 & 4.3056e-06 & 1.49  & 1.0349e-08 & 5.42 & 1.7664e-10 & 4.10 \\
\end{tabular}
\caption{\qq{\textit{(top)} Convergence on Manifold A of our GMLS solver based on the biharmonic formulation of the hydrodynamics in equation~\ref{equ_unsplit_stokes}.  The target sampling distance $h$ is discussed in Appendix~\ref{appendix:sampling_res_manifolds}.  \textit{(bottom)} Convergence on Manifold A of our GMLS solver based on the split formulation of the hydrodynamics in equation~\ref{equ_split_stokes}. }}
\label{table:hydro_conv_split_A}
\label{table:hydro_conv_biharmonic_A}
\label{table:hydro_conv_bh_and_split_A}
\end{center}
\end{table}

\begin{table}[H]
\begin{center}
Biharmonic Formulation: Manifold B\\
\begin{tabular}{ccccccc}
\textbf{}  & \multicolumn{2}{c}{\textbf{m = 4}}      & \multicolumn{2}{c}{\textbf{m = 6}}      & \multicolumn{2}{c}{\textbf{m = 8}}      \\\cline{2-7}
\textbf{h} & \textbf{$\ell_2$-error} & \textbf{Rate} & \textbf{$\ell_2$-error} & \textbf{Rate} & \textbf{$\ell_2$-error} & \textbf{Rate} \\
0.1    & 3.1890e-01 & -     & 1.0457e-01 & -    & 1.6845e+00 & -    \\
0.05   & 3.1951e-01 & -0.002 & 7.4388e-03 & 3.81 & 1.9954e-02 & 6.40 \\
0.025  & 2.4571e-02 & 3.69 & 1.2081e-03 & 2.62 & 2.9917e-04 & 6.05 \\
0.0125 & 5.6309e-03 & 2.12  & 6.9269e-05 & 4.11 & 2.6601e-05 & 3.48 \\
\end{tabular}
\\
\vspace{0.2cm}
Split Formulation: Manifold B\\
\begin{tabular}{ccccccc}
\textbf{}  & \multicolumn{2}{c}{\textbf{m = 4}}      & \multicolumn{2}{c}{\textbf{m = 6}}      & \multicolumn{2}{c}{\textbf{m = 8}}      \\\cline{2-7}
\textbf{h} & \textbf{$\ell_2$-error} & \textbf{Rate} & \textbf{$\ell_2$-error} & \textbf{Rate} & \textbf{$\ell_2$-error} & \textbf{Rate} \\
0.1    & 9.7895e-02 & -     & 6.5222e-02 & -    & 2.8024e-01 & -    \\
0.05   & 1.4383e-02 & 2.77 & 2.8402e-03 & 4.52 & 1.2100e-02 & 4.53 \\
0.025  & 3.6243e-03 & 1.98  & 3.9929e-04 & 2.82 & 4.9907e-04 & 4.59 \\
0.0125 & 7.8747e-04 & 2.20  & 1.2357e-05 & 5.00 & 5.7023e-06 & 6.44 \\
\end{tabular}
\caption{\qq{\textit{(top)} Convergence on Manifold B of our GMLS solver based on the biharmonic formulation of the hydrodynamics in equation~\ref{equ_unsplit_stokes}. \textit{(bottom)} Convergence on Manifold B of our GMLS solver based on the split formulation of the hydrodynamics in equation~\ref{equ_split_stokes}.}}
\label{table:hydro_conv_biharmonic_B}
\label{table:hydro_conv_split_B}
\label{table:hydro_conv_bh_and_split_B}
\end{center}
\end{table}

\begin{table}[H]
\label{Fountain_Stokes_Conv}
\begin{center}
Biharmonic Formulation: Manifold C\\
\begin{tabular}{ccccccc}
\textbf{}  & \multicolumn{2}{c}{\textbf{m = 4}}      & \multicolumn{2}{c}{\textbf{m = 6}}      & \multicolumn{2}{c}{\textbf{m = 8}}      \\\cline{2-7}
\textbf{h} & \textbf{$\ell_2$-error} & \textbf{Rate} & \textbf{$\ell_2$-error} & \textbf{Rate} & \textbf{$\ell_2$-error} & \textbf{Rate} \\
0.1    & 2.9886e+00 & -     & 8.0650e-01 & -    & 3.3799e-01 & -    \\
0.05   & 1.2926e+00 & 1.21 & 2.3277e-01 & 1.79 & 1.0993e+00 & -1.70 \\
0.025  & 2.8576e-01 & 2.18 & 2.1497e-02 & 3.44 & 7.1166e-03 & 7.28 \\
0.0125 & 4.2226e-02 & 2.76  & 1.4986e-03 & 3.84 & 9.8921e-05 & 6.17 \\
\end{tabular}
\\
\vspace{0.2cm}
Split Formulation: Manifold C\\
\begin{tabular}{ccccccc}
\textbf{}  & \multicolumn{2}{c}{\textbf{m = 4}}      & \multicolumn{2}{c}{\textbf{m = 6}}      & \multicolumn{2}{c}{\textbf{m = 8}}      \\\cline{2-7}
\textbf{h} & \textbf{$\ell_2$-error} & \textbf{Rate} & \textbf{$\ell_2$-error} & \textbf{Rate} & \textbf{$\ell_2$-error} & \textbf{Rate} \\
0.1    & 1.1346e+00 & -     & 8.8130e+01 & -    & 4.6473e+00 & -    \\
0.05   & 7.7801e-02 & 3.86 & 1.0276e-02 & 13.0 & 3.7375e-02 & 6.96 \\
0.025  & 1.6751e-02 & 2.22  & 1.8764e-03 & 2.45 & 4.2722e-04 & 6.46 \\
0.0125 & 1.7381e-03 & 3.27  & 4.2181e-05 & 5.48 & 9.1845e-06 & 5.54 \\
\end{tabular}
\caption{\qq{\textit{(top)} Convergence on Manifold C of our GMLS solver based on the biharmonic formulation of the hydrodynamics in equation~\ref{equ_unsplit_stokes}. \textit{(bottom)} Convergence on Manifold C of our GMLS solver based on the split formulation of the hydrodynamics in equation~\ref{equ_split_stokes}.}}
\label{table:hydro_conv_biharmonic_C}
\label{table:hydro_conv_split_C}
\label{table:hydro_conv_bh_and_split_C}
\end{center}
\end{table}

\begin{table}[H]
\label{Torus_Stokes_Conv}
\begin{center}
Biharmonic Formulation: Manifold D\\
\begin{tabular}{ccccccc}
\textbf{}  & \multicolumn{2}{c}{\textbf{m = 4}}      & \multicolumn{2}{c}{\textbf{m = 6}}      & \multicolumn{2}{c}{\textbf{m = 8}}      \\\cline{2-7}
\textbf{h} & \textbf{$\ell_2$-error} & \textbf{Rate} & \textbf{$\ell_2$-error} & \textbf{Rate} & \textbf{$\ell_2$-error} & \textbf{Rate} \\
0.08 & 3.3170e-01 & -     & 1.5154e-01 & -    & 1.2223e+01 & -    \\
0.04 & 2.4421e-02 & 3.82 & 4.6233e-03 & 5.11 & 3.9632e-03 & 11.7 \\
0.02 & 4.5705e-03 & 2.44 & 3.0246e-04 & 3.97 & 4.2784e-05 & 6.60 \\
0.01 & 1.4748e-03 & 1.62 & 1.9067e-05 & 3.96 & 5.4137e-07 & 6.26 \\
\end{tabular}
\\
\vspace{0.2cm}
Split Formulation: Manifold D\\
\begin{tabular}{ccccccc}
\textbf{}  & \multicolumn{2}{c}{\textbf{m = 4}}      & \multicolumn{2}{c}{\textbf{m = 6}}      & \multicolumn{2}{c}{\textbf{m = 8}}      \\\cline{2-7}
\textbf{h} & \textbf{$\ell_2$-error} & \textbf{Rate} & \textbf{$\ell_2$-error} & \textbf{Rate} & \textbf{$\ell_2$-error} & \textbf{Rate} \\
0.08 & 1.7719e-02 & -     & 1.4221e-02 & -    & 6.6061e+00 & -    \\
0.04 & 1.5473e-03 & 3.57 & 1.2632e-04 & 6.92 & 1.3431e-04 & 15.8 \\
0.02 & 1.3575e-04 & 3.54  & 3.2125e-06 & 5.35 & 5.0041e-07 & 8.15 \\
0.01 & 2.5891e-05 & 2.37  & 1.9018e-07 & 4.05 & 4.5906e-09 & 6.72 \\
\end{tabular}
\caption{\qq{\textit{(top)} Convergence on Manifold D of our GMLS solver based on the biharmonic formulation of the hydrodynamics in equation~\ref{equ_unsplit_stokes}. \textit{(bottom)} Convergence on Manifold D of our GMLS solver based on the split formulation of the hydrodynamics in equation~\ref{equ_split_stokes}.}}
\label{table:hydro_conv_biharmonic_D}
\label{table:hydro_conv_split_D}
\label{table:hydro_conv_bh_and_split_D}
\end{center}
\end{table}

We emphasize that these convergence studies take into account the full pipeline of our GMLS numerical methods as discussed in Section~\ref{sec:numerical_solver} and shown in Figure~\ref{fig:solver_schematic}.  This involves not only the solution of biharmonic or split equations, but also the GMLS reconstruction of the surface velocity field $\mb{v}$ from the computed vector-potential $\Phi$ and the calculation of the vector-potentials $\Psi = -\star\mb{d} \mb{b}$ for the body force density $\mb{b}$ which drives the flow.  These steps also each have a non-linear dependence on the geometry which contributes through our GMLS reconstructions from the point set sampling of the manifold as discussed in Section~\ref{sec:geo_reconstruct}.

In the convergence studies, we find in all cases that the GMLS solvers are able to resolve the surface hydrodynamic fields to a high level of precision.  The Manifolds $B$ and $C$ presented the most challenges for the solvers with largest prefactors in their convergence.  This is expected given the increased amount of resolution required to resolve the geometric contributions to the differential operators in the hydrodynamic equations~\ref{equ_unsplit_stokes}--~\ref{equ_split_stokes}.  In all cases, we found our GMLS solvers based on the split formulation performed better when using equation~\ref{equ_split_stokes} relative to our GMLS solvers based on the biharmonic formulation of equation~\ref{equ_unsplit_stokes}.  Interestingly, for Manifold $B$ and $C$ these differences for $m = 8$ where not as pronounced, see Table~\ref{table:hydro_conv_bh_and_split_B}--~\ref{table:hydro_conv_bh_and_split_C}.  We think this is a manifestation of the challenges in capturing the geometric contributions to the differential operator that with limited resolution will not benefit as much from the higher order approximations or split formulations relative to the case of less complicated geometries.  

We find in the case of Manifold $A$ that the GMLS solver for sufficiently large order ($m \geq 6$) converges at a rate of approximately $\sim 4^{th}$-order for the biharmonic formulation and at a rate of approximately $\sim 5^{th}$-order for the split formulation.  \qq{We base this on the overall trends, and some of this is a little obscured by the noise of the convergence after acheiving a high level of accuracy.  We suspect the last upward point of the error observed for $m = 8$ for the biharmonic formulation is likely a consequence of the conditioning of the linear system becoming a limiting factor.}  \qq{We note the overall high level of precision already achieved by that data point with errors on the order of $10^{-8}$, see Table~\ref{table:hydro_conv_bh_and_split_A}. We find there is a particular advantage of our GMLS solvers when based on the split formulation.  Our GMLS methods in this case are able to converge to much higher levels of precision achieving errors on the order $10^{-10}$ in the case of $m = 8$ at the largest resolutions considered, see Table~\ref{table:hydro_conv_bh_and_split_A}.}

\qq{In summary, our results show that both formulations of the GMLS solvers are able to achieve high-order convergence rates in approximating the hydrodynamic fields.  We emphasize that these results assess contributions from the entire pipeline that includes not only the GMLS solve but also the pre-processing and post-processing steps involving the curl operators that arise in our vector-potential formulation for incompressible hydrodynamic flows.}

\section{Conclusions}
\label{sec_conclusions}

We have developed high-order numerical methods for solving partial differential equations on manifolds.  Our apporach is based on GMLS approximations of the manifold shape, operators arising in differential geometry, and operators of differential equations.  We have introduced exterior calculus based approaches for generalizing the operators of vector calculus and techniques from mechanics to the context of manifolds.  Using this approach, we have formulated incompressible hydrodynamic equations for flows on curved surfaces.  We have also shown how our approaches in general can be used to formulate equations in terms of vector potentials facilitating development of other physical models with constraints to obtain numerical solvers.  We showed there are a few different ways to formulate vector-valued surface equations facilitating the development of GMLS solvers.  By comparisons with high precision manufactured solutions, we characterized our GMLS solvers and found they each exhibit high-order convergence rates in approximating manifold operators and in resolving hydrodynamics flows on surfaces.  We found the split formulations involving lower order differential operators to have particular advantages exhibiting the highest orders of convergence.  Many of our GMLS methods and exterior calculus approaches also can be utilized for the development of high-order solvers for other scalar-valued and vector-valued partial differential equations on manifolds.

\section*{Acknowledgments}
\label{sec_acknowledgements}
We acknowledge support to P.J.A. and B.J.G. from research grants DOE ASCR PhILMS DE-SC0019246 and NSF Grant DMS-1616353.  We also acknowledge the UCSB Center for Scientific Computing NSF MRSEC DMR-1121053. 

\bibliographystyle{plain}
\bibliography{paper_database}{}

\appendix

\section*{Appendix}
\label{sec_appendix}

\section{Operators on Manifolds, Monge-Gauge Parameterization, and Coordinate Expressions}
\label{appendix:monge_gauge}

To compute in practice the action of our operators during the GMLS reconstruction of the geometry of the manifolds or differential operators on scalar and vector fields on the surface, we use local Monge-Gauge parameterizations of the surface.  To obtain high-order accuracy we further expand expressions involving derivatives of the metric and other fields explicitly using symbolic algebra packages, such as Sympy~\cite{Sympy2017}.  This allows us to avoid some of the tedium notorous in differential geometry and to precompute offline the needed expressions for the action of our operators.  We summarize here the basic differential geometry of surfaces expressed in the Monge-Gauge and the
associated expressions we use in such calculations.

\subsection{Monge-Gauge Surface Parameterization}
\label{appendix:monge_gauge_param}
In the Monge-Gauge we parameterize locally a smooth surface in terms of the tangent plane coordinates $u,v$ and the height of the surface above this point as the function $h(u,v)$.  This gives the embedding map
\begin{eqnarray}
\label{equ:manifoldParam}
\mb{x}(u,v) = \bs{\sigma}(u, v) =  (u, v, \h(u,v)).
\end{eqnarray} 
We see that this parameterization of the surface is closely related to equation~\ref{eqn:chart}.  We can use the Monge-Gauge equation~\ref{equ:manifoldParam} to derive explicit expressions for geometric quantities.  The derivatives of $\bsy{\sigma}$ provide a basis $\partial_u, \partial_v$ for the tangent space as
\begin{eqnarray}
\label{eqn:sigmaU}
\partial_u = \bs{\sigma}_{u}(u, v) & = & (1, 0, \h_u(u, v) ) \\
\label{eqn:sigmaV}
\partial_v = \bs{\sigma}_{v}(u, v) & = & (0, 1, \h_v(u, v) ).
\end{eqnarray}
The first fundamental form $\mathbf{I}$ (metric tensor) and second fundamental form $\mathbf{II}$ (curvature tensor) are given by 
\begin{align}
\mathbf{I} = \begin{bmatrix}
E & F \\
F & G
\end{bmatrix} = \begin{bmatrix}
\bs{\sigma}_{u} \cdot \bs{\sigma}_{u} & \bs{\sigma}_{u} \cdot \bs{\sigma}_{v} \\
\bs{\sigma}_{v} \cdot \bs{\sigma}_{u} & \bs{\sigma}_{v} \cdot \bs{\sigma}_{v}
\end{bmatrix} =  \begin{bmatrix}
1 + \h_u(u,v)^2 & \h_u \h_v(u,v) \\
\h_u(u,v) \h_v(u,v) & 1 + \h_v(u,v)^2
\end{bmatrix}. \label{equ:I_mat}
\end{align}
and 
\begin{align}
\mathbf{II} = \begin{bmatrix}
L & M \\
M & N
\end{bmatrix} = \begin{bmatrix}
\bs{\sigma}_{u u} \cdot \bs{n} & \bs{\sigma}_{u v} \cdot \bs{n} \\
\bs{\sigma}_{v u} \cdot \bs{n} & \bs{\sigma}_{v v}\cdot \bs{n}
\end{bmatrix}= \frac{1}{\sqrt{ 1 + \h_u^2 + \h_v^2 }}\begin{bmatrix}
\h_{uu} & \h_{uv} \\
\h_{uv} & \h_{vv}
\end{bmatrix} . \label{equ:II_mat} 
\end{align}
The $\mb{n}$ denotes the outward normal on the surface and is given by
\begin{eqnarray}
\bs{n}(u, v) = \frac{\bs{\sigma}_{u}(u, v) \times \bs{\sigma}_{v}(u, v)}{\| \bs{\sigma}_{u}(u, v) \times \bs{\sigma}_{v}(u, v) \|} = \frac{1}{\sqrt{1+\h_u^2 + \h_v^2}}(-\h_u, -\h_v, 1). \label{equ:normal_vec_def}
\end{eqnarray}
We use throughout the notation for the metric tensor $\mb{g} = \mathbf{I}$  interchangeably.  For notational convenience, we use the tensor notation for the metric tensor $g_{ij}$ and for its inverse $g^{ij}$.  These correspond to the first and second fundamental forms as
\begin{eqnarray}
g_{ij} = \left[\mathbf{I}\right]_{i,j}, \hspace{0.3cm} g^{ij} = \left[ \mathbf{I}^{-1}\right]_{i,j}.
\end{eqnarray}
For the metric tensor $\mb{g}$, we also use the notation $|g| = \det(\mb{g})$ and 
have that
\begin{eqnarray}
\sqrt{|g|} = \sqrt{\det(\mathbf{I})} = \sqrt{1 + \h_u^2 + \h_v^2} = \| \vec{\sigma}_{u}(u, v) \times \vec{\sigma}_{v}(u, v) \|.
\label{met_fac_def}
\end{eqnarray}
The provides the local area element as $dA_{u,v} = \sqrt{|g|}du dv$.  To compute quantities associated with curvature of the manifold we construct the Weingarten map~\cite{Pressley2001} which can be expressed as
\begin{eqnarray}
\mb{W} = \mb{I}^{-1} \mb{II}.
\end{eqnarray}
The Gaussian curvature $K$ can be expressed in the Monge-Gauge as
\begin{eqnarray}
K(u,v) = \det\left(\mb{W}(u,v)\right) = \frac{\h_{uu}\h_{vv} - \h_{uv}^2}{(1 + \h_u^2 + \h_v^2)^2}.
\end{eqnarray}
For further discussions of these tensors and more generally the differential geometry of manifolds see~\cite{Pressley2001,Abraham1988,SpivakDiffGeo1999}.  We use these expressions as the basis of our calculations of the action of our surface operators.

\subsection{Coordinate Expressions for Surface Operators}
\label{appendix:monge_gauge_coord_exp}
We use local Monge-Gauge parameterizations of the manifold to compute the geometric operators needed in our surface hydrodynamic equations.  Consider the negative semi-definite scalar Laplace-Beltrami operator that acts on $0$-forms which can be expressed as $\Delta_{LB} = -\bs{\delta} \mb{d} = -\Delta_H$, where $\Delta_H$ is the Hodge Laplacian.  This operator can be expressed in coordinates as
\begin{equation}
\Delta_{LB} = \frac{1}{\sqrt{|g|}}\partial_i \left(g^{ij}\sqrt{|g|} \partial_j \right).
\end{equation}
The $g_{ij}$ denotes the metric tensor, $g^{ij}$ the inverse metric tensor, and $|g|$ the determinant of the metric tensor as in Appendix~\ref{appendix:monge_gauge_param}.  For the Monge-Gauge parameterization $(u,v)$, we find it useful to consider
\begin{eqnarray}
\label{eqn:ell_grouping}
\ell_{ij} = \left(\sqrt{|g|}g^{ij}\right)\partial_{ij} + \left(\partial_{i}\sqrt{|g|}g^{ij}\right)\partial_j.
\end{eqnarray}
We use the convention that $\partial_1 = \partial_u$ and $\partial_2 = \partial_v$.  This allows us to express 
\begin{eqnarray}
\Delta_{LB} = (1/\sqrt{|g|}) \sum_{ij} \ell_{ij}.
\end{eqnarray}  
We can further express the prefactor terms involving the metric appearing in equation~\ref{eqn:ell_grouping} as 
\begin{equation}
\label{equ_g_gij}
    \sqrt{|g|}g^{ij} = 
    \begin{cases}
       {g_{v v}}/{\sqrt{|g|}} = \frac{1+ \h_v^2}{\sqrt{1 + \h_u^2 + \h_v^2}} & \text{if: $i = j = u$} \\
       {g_{u u}}/{\sqrt{|g|}} = \frac{1+ \h_u^2}{\sqrt{1 + \h_u^2 + \h_v^2}} & \text{if: $i = j = v$} \\
       {-g_{u v}}/{\sqrt{|g|}} = {-g_{v u}}/{\sqrt{|g|}} = \frac{-\h_u \h_v}{\sqrt{1 + \h_u^2 + \h_v^2}} & \text{if: $i \ne j$}.
    \end{cases}
\end{equation}
The utility of these decompositions and expressions is that we can construct operators for GMLS approximation while avoiding the need to compose numerical differentiation procedures. This allows us to compute directly the action on the reconstruction space functions $p \in \mathbb{V}_h$.  This decomposition is also useful to help simplify symbolic expansions when we compute the Bi-Laplace-Beltrami operator $\Delta_{LB}^2$, which poses the most significant computational challenges in our current numerical calculations.  We compute in practice the Bi-Laplace-Beltrami operator $\Delta_{LB}^2$ using symbolic algebra system.

\subsection{Exterior Calculus Operators Expressed in Coordinates}
\label{appendix:monge_gauge_coord_ext_calc_exp}
In our notations throughout, we take the conventions that for differential $0$-forms (scalar functions) $f_j = \partial_{x^j} f$, for differential $1$-forms (co-vector fields) $\bs{\alpha} = \alpha_j \mb{d}{x}^j$, and for vector fields $\mb{v}=v^{j}\partial_{j}$.  In each case we have $j \in \{u,v\}$.  The isomorphisms $\sharp$ and $\flat$ between vectors and co-vectors can be expressed explicitly as
\begin{eqnarray}
\mb{v}^{\flat} &=& (v^{u}\bs{\sigma}_{u} + v^{v}\bs{\sigma}_{v})^{\flat}\\
\nonumber
 &=& v^{u}g_{u u}\mb{d}u + v^{u}g_{u v}\mb{d}v + v^{v}g_{v u}\mb{d}u + v^{v}g_{v v}\mb{d}v \\
 \nonumber
  &=& (v^{u}g_{u u} + v^{v}g_{v u})\mb{d}u + (v^{u}g_{u v} + v^{v}g_{v v})\mb{d}v \\ 
\nonumber
\\
(\bs{\alpha})^{\sharp} &=& (\alpha_{u}\mb{d}u + \alpha_{v}\mb{d}v)^{\sharp} \\
\nonumber
 &=& \alpha_{u}g^{u u}\bs{\sigma}_{u} + \alpha_{u}g^{u v}\bs{\sigma}_{v} + \alpha_{v}g^{v u}\bs{\sigma}_{u} + \alpha_{v}g^{v v}\bs{\sigma}_{v} \\
 \nonumber
 &=& (\alpha_{u}g^{u u} + \alpha_{v}g^{v u})\bs{\sigma}_{u} + (\alpha_{u}g^{u v} + \alpha_{v}g^{v v})\bs{\sigma}_{v}
\end{eqnarray}
We use the conventions for the notation that for the embedding map $\bs{\sigma}$ we have $\bs{\sigma}_u = \partial_{u}$ and $\bs{\sigma}_v = \partial_{v}$ as in Appendix~\ref{appendix:monge_gauge_param}.  The exterior derivatives on these $k$-forms can be expressed as
\begin{eqnarray}
\mb{d}f &=& (\partial_{u}f) \mb{d}u + (\partial_{v} f) \mb{d}v = f_{u}\mb{d}u + f_{v}\mb{d}v \\
\mb{d}\bsy{\alpha} &=& (\partial_{u} \alpha_{v} - \partial_{v} \alpha_{u}) \mb{d}u \wedge \mb{d}v.
\end{eqnarray}
The generalized curl of a $0$-form and $1$-form can be expressed in coordinates as
\begin{eqnarray}
-\star \mb{d}f & = & 
\sqrt{|g|}(f_{u}g^{u v} + f_{ v} g^{v v})\mb{d}u - \sqrt{|g|}(f_{u}g^{u u} + f_{v}g^{v u})\mb{d}v \\
-\star \mb{d} \alpha & = & \frac{\partial_{v} \alpha_{u} - \partial_{u} \alpha_{v} }{\sqrt{|g|}}.
\end{eqnarray}
Combining the above equations, we can express the generalized curl as
\begin{eqnarray}
\label{equ_gen_curl_0form}
(-\star \mb{d}f)^{\sharp} &= & \mbox{curl}_\mathcal{M}(f)\\ \nonumber 
&=&   ([\sqrt{|g|}(f_{u}g^{u v} + f_{ v} g^{v v})]g^{u u} + [- \sqrt{|g|}(f_{u}g^{u u} + f_{v}g^{v u})]g^{v u})\bs{\sigma}_{u}\\ 
\nonumber
&+ & ([\sqrt{|g|}(f_{u}g^{u v} + f_{ v} g^{v v})]g^{u v} + [- \sqrt{|g|}(f_{u}g^{u u} + f_{v}g^{v u})]g^{v v})\bs{\sigma}_{v} \\
\nonumber
& = & \frac{f_{v}}{\sqrt{|g|}} \bs{\sigma}_{u} - \frac{f_{u}}{\sqrt{|g|}} \bs{\sigma}_{v} \\ 
\nonumber
\\
\label{equ_gen_curl_1form}
-\star \mb{d}\mb{v}^{\flat} & = & \mbox{curl}_\mathcal{M}(\mb{v}) = \frac{\partial_{v}(v^{u}g_{u u} + v^{v}g_{v u}) - \partial_{u} (v^{u}g_{u v} + v^{v}g_{v v}) }{\sqrt{|g|}}.
\end{eqnarray}

We also mention that the velocity field of the hydrodynamic flows $\mb{v}$ is recovered from the vector potential $\Phi$ as $\mb{v}^{\flat} = -\star \mb{d} \Phi$.  We obtain the velocity field as $\mb{v} = \left(\mb{v}^{\flat} \right)^{\sharp} = \left(-\star \mb{d} \Phi\right)^{\sharp}$ using equation~\ref{equ_gen_curl_0form}.  Similarly from the force density $\mb{b}$ acting on the fluid, we obtain from equation~\ref{equ_gen_curl_1form} the vector potential for the force density as $\Psi = -\star\mb{d}\mb{b}^{\flat}$.  This is used in the vector-potential formulation of the hydrodynamics in equation~\ref{equ_unsplit_stokes} and equation~\ref{equ_split_stokes}.  We expand these expressions further as needed in coordinates using symbolic algebra methods.   This provides the needed expressions for computing these operations.  Additional details and discussions of these operators and our overall approach also can be found in our related papers~\cite{AtzbergerSoftMatter2016,AtzbergerGrossSurfExtCalc2017}.

\section{Convergence Results for the Generalized Curl Operators}
\label{appendix:conv_curl}

\qq{We report tabulated results for the GMLS approximations of the operators $\mathcal{L}_{C0}$ and $\mathcal{L}_{C1}$ discussed in Section~\ref{sec:conv_ops}.}

\begin{table}[H]
\label{C0_Conv}
\begin{center}
\begin{tabular}{ccccccccccc}
\textbf{}  & \multicolumn{2}{c}{\textbf{Manifold A}}      & \multicolumn{2}{c}{\textbf{Manifold B}} & \multicolumn{2}{c}{\textbf{Manifold C}} & & & \multicolumn{2}{c}{\textbf{Manifold D}}      \\\cline{2-7}\cline{10-11}
\textbf{h} & \textbf{$\ell_2$-error} & \textbf{Rate} & \textbf{$\ell_2$-error} & \textbf{Rate} & \textbf{$\ell_2$-error} & \textbf{Rate} & & \textbf{h} & \textbf{$\ell_2$-error} & \textbf{Rate} \\
0.1 & 2.7152e-05 & - & 1.5075e-03 & -  & 4.8243e-01  & -&  & .08 &  2.1570e-03 & -  \\
0.05 & 3.8309e-07 & 6.07 & 3.0281e-05 & 5.64 & 2.4465e-04 & 10.9 & & .04 & 2.2565e-05 & 6.68 \\
0.025 & 5.8491e-09 & 6.00 & 6.9649e-07 & 5.43 & 6.1779e-06 & 5.31 & & .02 &  3.3550e-07 & 6.13 \\
0.0125 & 8.8291e-11 & 6.04 & 1.3078e-08 & 5.72 & 1.1817e-07 &  5.71 & & .01 & 4.9708e-09 & 6.04 \\
\end{tabular}
\end{center}
\caption{Convergence of GMLS Approximation of the Surface Curl Operator on Scalars $\mathcal{L}_{C0}$.  We use GMLS with ($k = 1, m = 6)$ and find the methods have $\sim 6^{th}$-order asymptotic convergence.}
\end{table}

\begin{table}[H]
\label{C1_Conv}
\begin{center}
\begin{tabular}{ccccccccccc}
\textbf{}  & \multicolumn{2}{c}{\textbf{Manifold A}}      & \multicolumn{2}{c}{\textbf{Manifold B}} & \multicolumn{2}{c}{\textbf{Manifold C}} & & & \multicolumn{2}{c}{\textbf{Manifold D}}      \\\cline{2-7}\cline{10-11}
\textbf{h} & \textbf{$\ell_2$-error} & \textbf{Rate} & \textbf{$\ell_2$-error} & \textbf{Rate} & \textbf{$\ell_2$-error} & \textbf{Rate} & & \textbf{h} & \textbf{$\ell_2$-error} & \textbf{Rate} \\
0.1 & 9.2312e-04 & - & 1.5887e-02 & -  & 5.2497e+01  & -&  & .08 &  1.9686e-02 & -  \\
0.05 & 1.4851e-05 & 5.88 & 1.2736e-03 & 3.64 & 1.3126e-02 & 8.65 & & .04 & 2.0410e-04 & 6.70 \\
0.025 & 2.3374e-07 & 5.96 & 1.2597e-04 & 3.33 & 5.6087e-04 & 4.55 & & .02 &  3.0223e-06 & 6.13 \\
0.0125 & 3.5970e-09 & 6.01 & 5.1267e-06 & 4.61 & 1.4082e-05 &  5.32 & & .01 & 4.3847e-08 & 6.07 \\
\end{tabular}
\end{center}
\caption{Convergence of GMLS Approximation of the Surface Curl Operator on Vectors $\mathcal{L}_{C1}$.  We use GMLS with ($k_1 = 2, k_2 = 1, m = 6)$ and find the methods have $\sim 5^{th}$-order asymptotic convergence or greater.  It is notable that in the case of Manifold $A$ and $D$ we in fact see $\sim 6^{th}$-order convergence.  This manifests since the manifolds have a relatively symmetric geometry compared to Manifold $B$ and $C$, see Figure~\ref{fig:manifolds}.  This results in a simplification with fewer non-zero terms and derivatives associated with the contributions of the geometry to the operator.  As a consequence, the GMLS approximation at a given order $m$ becomes more accurate by one order for Manifold $A$ and $D$.
}
\end{table}

The Manifolds B and C have more complicated geometry and require more resolution to see behaviors in the asymptotic regime with a high-degree basis.  We see that by lowering the degree of the basis these operators exhibit more readily behaviors in the asymptotic regime in Table~\ref{table:conv_C0_lower} and~\ref{table:conv_C1_lower}.

\begin{table}[H]
\label{C0_Conv_lower_deg}
\begin{center}
\begin{tabular}{ccccccc}
\textbf{}  & \multicolumn{2}{c}{\textbf{Manifold B}} & \multicolumn{2}{c}{\textbf{Manifold C}} \\\cline{2-5}
\textbf{h} & \textbf{$\ell_2$-error} & \textbf{Rate} & \textbf{$\ell_2$-error} & \textbf{Rate} \\
0.1 &  5.2558e-03 & -  & 1.2083e-02  & - \\
0.05 &  3.6359e-04 & 3.85 & 1.0345e-03 & 3.54 \\
0.025 &  2.3078e-05 & 3.97 & 7.3790e-05 & 3.81  \\
0.0125 & 1.4569e-06 & 3.98 & 4.8316e-06 &  3.93  \\
\end{tabular}
\end{center}
\caption{Convergence of GMLS Approximation of the Surface Curl on Scalars $\mathcal{L}_{C0}$.  We use GMLS with ($k = 1, m = 4)$ and find the methods have $\sim 4^{th}$-order asymptotic convergence.}
\label{table:conv_C0_lower}
\end{table}

\begin{table}[H]
\label{C1_Conv_lower_deg}
\begin{center}
\begin{tabular}{ccccccc}
\textbf{}  & \multicolumn{2}{c}{\textbf{Manifold B}} & \multicolumn{2}{c}{\textbf{Manifold C}} \\\cline{2-5}
\textbf{h} & \textbf{$\ell_2$-error} & \textbf{Rate} & \textbf{$\ell_2$-error} & \textbf{Rate} \\
0.1 &  6.3586e-01 & -  & 7.6579e-01  & - \\
0.05 &  1.6568e-01 & 1.94 & 2.1680e-01 & 1.82 \\
0.025 & 4.1633e-02 & 1.99 & 5.6498e-02 & 1.94  \\
0.0125 & 1.0399e-02 & 1.99 & 1.4336e-02 &  1.98  \\
\end{tabular}
\end{center}
\caption{Convergence of GMLS Approximation of the Surface Curl on Vectors $\mathcal{L}_{C1}$.  We use GMLS with ($k = 1, m = 2)$ and find the methods have $\sim 2^{nd}$-order asymptotic convergence.}
\label{table:conv_C1_lower}
\end{table}

\section{Dependence of GMLS Approximations on the Point Sampling}
\label{appendix:stability_sampling}

\begin{figure}[H]  
\centerline{\includegraphics[width=0.6\columnwidth]{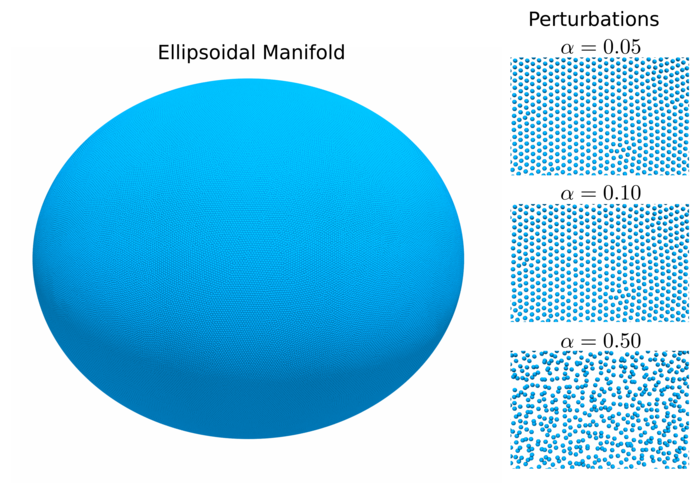}}
\caption{\qq{Ellipsoidal Manifold, Point Samples, and Perturbations.   \textit{(left)} The ellipsoidal manifold with $(x/\ell_x)^2 + (y/\ell_y)^2 + (z/\ell_z)^2  = 1$ with $\ell_x = 1.2$, $\ell_y = 1.2$, $\ell_z = 1$.  The manifold is sampled with $154,182$ points.  \textit{(right)} The points are perturbed by noise having strength $\alpha = 0.05,0.10,0.50$.}
}
\label{fig:stability_study}	
\end{figure}
\qq{We investigate the dependence of the GMLS approximations on the point sampling.  We perform studies of the sampling of the manifold perturbed by noise.  As discussed in  Appendix~\ref{appendix:sampling_res_manifolds}, we start with the manifold sampled by a collection of points $\{x_i\}$ using DistMesh~\cite{DistMesh2004}.  These are nearly uniform as characterized by equation~\ref{eqn:def_q_uniform}.  We then add Gaussian noise $\eta_i$ with mean zero and standard derivation of $\sigma = \alpha \ell_*$ to each of the points $x_i$.  The $\ell_*$ is $1/3$ the smallest nearest neighbor distance.
We project these points back to the manifold to obtain the new perturbed sampling $\tilde{x}_i = \mathcal{P}(x_i + \eta_i)$, where $\mathcal{P}$ denotes the projection mapping.  For the radial manifold shapes, $\mathcal{P}$ is the mapping back to the manifold along the radial directions.  We study the perturbations when $\alpha \in \{0.05, 0.1, 0.5\}$ for the ellipsoid $(x/\ell_x)^2 + (y/\ell_y)^2 + (z/\ell_z)^2  = 1$ with $\ell_x = 1.2$, $\ell_y = 1.2$, $\ell_z = 1$.  We show realizations of the point sampling and perturbations in Figure~\ref{fig:stability_study} and the results of our convergence studies in Table~\ref{table:stability_study}.
}

\qq{
\begin{table}[H]
\begin{center}
\begin{tabular}{ccccccccc}
\textbf{}  & \multicolumn{2}{c}{\textbf{Laplace-Beltrami}} & \multicolumn{2}{c}{\textbf{Biharmonic}} \\\cline{2-5}
\textbf{$\alpha$} & \textbf{$\ell_2$-error} & \textbf{Rate} & \textbf{$\ell_2$-error} & \textbf{Rate} \\
0.00 & 3.437-04 & 2.03 & 2.828-04 & 2.54 \\
0.05 & 3.425-04 & 2.03 & 8.655-04 & 2.20 \\
0.10 & 3.419-04 & 2.04 & 9.032-04 & 2.17 \\
0.50 & 3.650-04 & 2.07 & 1.464-03 & 2.18  \\
\end{tabular}
\end{center}
\caption{\qq{Dependence of GMLS Approximations on the Point Sampling.  We investigate the convergence when perturbing the manifold sample points by $\tilde{x}_i = \mathcal{P}(x_i + \eta_i)$ where the Gaussian noise $\eta_i$ has mean zero and standard deviation of $\sigma = \alpha \ell_*$.  The $\ell_*$ is $1/3$ of the smallest nearest neighbor distance for the unperturbed sampling.  We consider the ellipsoidal manifold with $(x/\ell_x)^2 + (y/\ell_y)^2 + (z/\ell_z)^2  = 1$ with $\ell_x = 1.2$, $\ell_y = 1.2$, $\ell_z = 1$ sampled with the number of sample points $2,350$, $9,566$, $38,486$, and $154,182$.  We study the accuracy of the solvers for $\mathcal{L}u = -f$ where $\mathcal{L}$ is the Laplace-Beltrami Operator and the Biharmonic Operator.  The solver for the Biharmonic Operator uses the split formulation.  The function $f$ is generated using angular coordinates to obtain the real-part of the spherical harmonic $Y_4^5$ projected to the ellipsoidal surface.  The ellipsoid considered is given by $(x/\ell_x)^2 + (y/\ell_y)^2 + (z/\ell_z)^2  = 1$ with $\ell_x = 1.2$, $\ell_y = 1.2$, $\ell_z = 1$.}}
\label{table:stability_study}
\end{table}
}

\qq{From Table~\ref{table:stability_study}, we see that the GMLS methods are robust to perturbations in the sampling both for the Laplace-Beltrami Operator and for the Biharmonic Operator.  We remark that these solvers each have a non-trivial dependence not only on the approximation of the differential operators but also in the approximations performed for the local geometry of the manifold.  This indicates in GMLS that the local least-squares fitting is not overly sensitive to the point sampling that is used in constructing the approximations.}

\section{Sampling Resolution of the Manifolds}
\label{appendix:sampling_res_manifolds}
\qq{A summary of the sampling resolution $h$ used for each of the manifolds is provided in Table~\ref{table:sampling_manifolds}.  We refer to $h$ as the \textit{target fill distance}.  For each of the manifolds, we achieve a nearly uniform collection of the points as in equation~\ref{eqn:def_q_uniform} using DistMesh~\cite{DistMesh2004}.}   We emphasize this approach was used only for convenience to obtain  quasi-uniform samplings and other sampling techniques can also be utilized for this purpose of representing the manifolds.  We specify $h$ and the algorithm produces a point sampling of the manifold.  In practice,  we have found this yields a point spacing with neighbor distances varying by only $\approx \pm 30 \%$ relative to the target distance $h$.  We summarize for each of the manifolds how this relates to the number of sample points $n$ in Table~\ref{table:sampling_manifolds}.  
\begin{table}[H]
\begin{center}
 \begin{tabular}{||c | c c |c c |c c |c c ||} 
 \hline
Refinement Level & \textbf{A}: h \hspace{0.22cm} & n & \textbf{B}: h \hspace{0.22cm} & n & \textbf{C}: h \hspace{0.22cm} & n & \textbf{D}: h \hspace{0.22cm} & n \\ 
 \hline\hline
 $1$ & .1 & 2350 & .1 & 2306 & .1 & 2002 & .08 & 1912 \\ 
 \hline
 $2$ & .05 & 9566 & .05 & 9206 & .05 & 7998 & .04 & 7478 \\
 \hline
 $3$ & .025 & 38486 & .025 & 36854 & .025 & 31898 & .02 & 29494 \\
 \hline
 $4$ & .0125 & 154182 & .0125 & 147634 & .0125 & 127346 & .01 & 118942 \\
 \hline
\end{tabular}
\captionof{table}{Sampling Resolution for each of the Manifolds A--D.  Relation between the target distance $h$ and the number of sample points $n$ used for each of the manifolds.  In each case, the neighbor distances between the points sampled were within $\approx \pm 30\%$ of the target distance $h$.  } 
\label{table:sampling_manifolds}
\end{center}
\end{table}

\end{document}